\documentclass[a4paper,leqno]{amsart}
\usepackage{amsfonts}
\theoremstyle{plain}
\begingroup

\newtheorem{theorem}{Theorem}[section]
\newtheorem{lemma}[theorem]{Lemma}
\newtheorem{proposition}[theorem]{Proposition}
\newtheorem{corollary}[theorem]{Corollary}
\newtheorem{remark}[theorem]{Remark}
\newtheorem{definition}[theorem]{Definition}
\endgroup

\theoremstyle{definition}
\theoremstyle{remark}

\mathchardef\emptyset="001F

\numberwithin{equation}{section}

\newcommand{\e}{\varepsilon}
\newcommand{\Om}{\Omega}
\newcommand{\Omk}{{\Om\setmeno K}}
\newcommand{\om}{\omega}
\newcommand{\R}{{\mathbb R}}
\newcommand{\wto}{{\rightharpoonup}}
\newcommand{\setmeno}{\!\setminus\!}

\newcommand{\huno}{{\mathcal H}^{1}}
\newcommand{\E}{{\mathcal E}}
\newcommand{\gdot}{{\dot{g}}}
\newcommand{\K}{{\mathcal K}(\overline\Om)}
\newcommand{\Kmf}{{\mathcal K}_m(\overline\Om)}
\newcommand{\Kmfbul}{{\mathcal K}'_m(\overline\Om)}
\newcommand{\KA}{{\mathcal K}_{q+m-h}(\overline A)}
\newcommand{\salt}{\noalign{\vskip 3pt}}

\title[Quasi-static growth of brittle fractures]
{A  model for the quasi-static growth\\
of brittle fractures based on local minimization
     }
\author[Gianni Dal Maso]{Gianni Dal Maso}
\address[Gianni Dal Maso]{SISSA, Via Beirut 2-4, 34014 Trieste,
Italy}
\email[Gianni Dal Maso]{dalmaso@sissa.it}
\author[Rodica Toader]{Rodica Toader}
\address[Rodica Toader]{Dip.to di Scienze Matematiche, Via A. Valerio 12/1,
34127, Trieste,
Italy}
\email[Rodica Toader]{toader@mathsun1.univ.trieste.it}

\begin{document}

\begin{abstract}
We study a variant of the variational model for the quasi-static 
growth of brittle fractures proposed by Francfort and Marigo in 
\cite{FraMar3}. The main feature of our model is that, in the 
discrete-time formulation, in each step we do not consider 
absolute minimizers of the energy, but, in a sense, 
we look for local minimizers which are 
sufficiently close to the approximate solution obtained in the 
previous step. This is done by
introducing in the variational problem an additional term which
penalizes the $L^2$-distance between the approximate solutions at two 
consecutive times. We study the continuous-time version of this 
model, obtained by passing to the limit as the time step tends to 
zero, and show that it satisfies (for almost every time) some minimality 
conditions which are slightly different from those considered in 
\cite{FraMar3} and \cite{DM-T}, but are still enough to prove
(under suitable regularity assumptions on the crack path)
that the classical Griffith's criterion holds at the crack tips. 
We prove also that, if no initial crack is present
and if the data of the problem are sufficiently smooth, 
no crack will develop 
in this model, provided the
penalization term is large enough.
\end{abstract}
\maketitle
{\small

\bigskip

\bigskip
}
\bigskip
\bigskip

\begin{section}{INTRODUCTION}

In this paper we present a variational model for the 
irreversible quasi-static growth of 
brittle fractures in the two-dimensional
antiplane case, subject to a time dependent
boundary displacement. 
The reference configuration is a bounded Lipschitz domain $\Omega$ 
of the plane, whose boundary $\partial\Om$ is divided 
into two disjoint locally connected subsets
$\partial_D\Om$ and  $\partial_N\Om$, where 
we prescribe a nonhomogeneous Dirichlet condition and a homogeneous
Neumann condition, 
respectively. 
According to Griffith's theory, the energy considered in the model is 
given by
\begin{equation}\label{energy}
E(u,K):=\int_{\Om\setminus K}|\nabla u|^2\,dx + \huno(K)\,,
\end{equation}
where the compact set $K\subset\overline\Om$
represents the crack in the reference 
configuration, the scalar function $u$ represents the displacement 
orthogonal to the plane of 
$\Om$, and $\huno$ is the one-dimensional Hausdorff measure.

For technical reasons, due to the behaviour of the solutions of 
Neumann problems in domains with cracks (see \cite{BucVar1}, 
\cite{DM-T}), we impose an a priori bound on the 
number of connected components of $K$.

Given a time dependent energy functional ${\mathcal F}(z,t)$, defined for $z$ 
in a Banach space $Z$ and for $t\in [0,T]$,
a quasi-static evolution corresponding to ${\mathcal F}$ 
is a function $t\mapsto z(t)$ which satisfies 
the equality $\nabla_z{\mathcal F}(z(t),t)=0$ for every $t\in [0,T]$.
A standard way to obtain this function is by singular perturbation. 
We 
consider the $\e$-gradient 
flow
\begin{equation}\label{gradflow}
\e\dot{z}^\e(t)+\nabla_z{\mathcal F}(z^\e(t),t)=0\,,
\end{equation}
starting from a local minimizer $z_0$ of ${\mathcal F}(\cdot,0)$.
Under suitable assumptions on ${\mathcal F}$ the solutions $z^\e(t)$ 
converge, as $\e\to 0$, to a function $z(t)$ which satisfies the equation
\begin{equation}\label{qs}
\nabla_z{\mathcal F}(z(t),t)=0\,;
\end{equation}
moreover, due to the choice of the sign in (\ref{gradflow}), it turns 
out that $z(t)$ is a local minimizer of ${\mathcal F}(\cdot,t)$ for a 
generic $t\in[0,T]$. 

Heuristically, the potential well of ${\mathcal 
F}(\cdot,0)$ corresponding to $z_0$ will be slightly deformed for $t$ 
small, and
the solution $z(t)$ of (\ref{qs}) 
obtained by this approximation method
follows the local minimizer of the deformed potential 
well. It may happen that for some critical value $t_0$ this potential 
well disappears, and for this special time $z(t_0)$ will be only 
a critical point of 
${\mathcal F}(\cdot,t_0)$; in general, in this case $z(t)$  is 
discontinuous at $t_0$ and jumps to another point $z(t_0+)$, 
which is a local minimizer of
${\mathcal F}(\cdot,t_0)$; the evolution  continues then  in
this new potential well. By a simple rescaling argument we see that
$z(t_0+)$ can also be obtained from $z(t_0)$ by solving 
the gradient flow
(\ref{gradflow}) with 
$\e=1$ and with initial conditions close to $z(t_0)$, and taking then the limit 
as $t\to+\infty$.

We want to adapt these ideas to the case of the energy (\ref{energy}) 
with the time dependent Dirichlet boundary condition $u(t)=g(t)$ on 
${\partial_D\Omk(t)}$ and with initial condition $(u_0, K_0)$. 
We assume that $g(t)$ is the trace on $\partial_D\Om$ of a function 
of $H^1(\Om)$, still denoted by $g(t)$, and that the map $t\mapsto 
g(t)$ belongs to $AC([0,T];H^1(\Om))\cap L^\infty([0,T];L^\infty(\Om))$.

Since we are
looking for equilibria, it is natural to assume that $u_0$ minimizes
$$
\int_{\Om\setminus K_0}|\nabla u|^2\,dx 
$$
among all functions $u\in H^1(\Omk_0)$ with $u=g(0)$ on 
$\partial_D\Omk_0$.

The main difficulty  in the definition of the $\e$-gradient flow 
for  (\ref{energy}) is that this energy
it is neither 
differentiable nor convex, due to the presence of the term 
$\huno(K)$, and therefore we can not rely on a notion of 
(sub)differential.
Following \cite{AmBr} and \cite{ChD}, we define the 
$\e$-gradient flow of the energy  (\ref{energy}) using an
approximation by a 
discrete-time process based on an implicit scheme. 

Let us fix an 
integer $m\ge 1$ and let $\Kmf$ be the set of all compact subsets $K$ 
of $\overline \Om$ with at most $m$ connected components and with 
$\huno(K)<+\infty$. We consider also the set $\Kmfbul$ of all $K\in 
\Kmf$ without isolated points, and we assume that the initial crack 
$K_0$ belongs to $\Kmfbul$.

Given the time step $\delta>0$, for any integer $i\ge 0$ let 
$t_i^\delta:=i\delta$, and, for $t_i^\delta\le T$, let $g_i^\delta:=g(t_i^\delta)$.
We define $(u_i^{\e,\delta}, K_i^{\e,\delta})$ inductively as follows:
$(u_{0}^{\e,\delta}, K_{0}^{\e,\delta}):=(u_0, K_0)$;
for $i\geq1$ we define $(u_i^{\e,\delta}, K_i^{\e,\delta})$ as a solution of 
the minimum problem
\begin{equation}\label{piepsdelta}
\min_{(u,K)}\big\{ E(u,K)+
\frac{\e}{\delta} \,\|u- u_{i-1}^{\e,\delta}\|^2
\big\}\,,
\end{equation}
where $\|\cdot\|$ denotes the $L^2$-norm in $\Om$, and the minimum is 
taken over all pairs $(u,K)$ such that
$K\in\Kmfbul$, $K\supset K_{i-1}^{\e,\delta}$, $u\in H^1(\Omk)$, 
$u=g_i^\delta$ on $\partial_D\Omk$. The constraint $K\supset 
K_{i-1}^{\e,\delta}$ reflects the irreversibility of the fracture 
process.

We define now the step functions $u^{\e,\delta}(t)$ and 
$K^{\e,\delta}(t)$
on $[0,T]$
by setting  $u^{\e,\delta}(t):=u_{i}^{\e,\delta}$ and 
$K^{\e,\delta}(t):=K_{i}^{\e,\delta}$,
 for
$t_i^\delta\leq t<t_{i+1}^\delta$.

The limit $(u^\e(t),K^\e(t))$ of $(u^{\e,\delta}(t),K^{\e,\delta}(t))$ along 
a suitable sequence 
$\delta_k\to0$ is by definition the $\e$-gradient flow 
for the energy (\ref{energy}). 
In order to obtain the quasi-static evolution for this energy by the 
singular perturbation approach, we should consider now the limit 
of $(u^\e(t),K^\e(t))$ as $\e\to0$ along a suitable sequence. 
This can be done, but we are not able to prove satisfactory properties 
of the limit evolution process $(u(t),K(t))$. 

Therefore we prefer to 
consider a variant of the singular perturbation method. We study the 
limit of $(u^{\e,\delta}(t),K^{\e,\delta}(t))$ as $\e$ and $\delta$ 
tend to zero simultaneously, with $\e$ proportional to~$\delta$. In 
particular, given $\lambda>0$, we assume that the coefficient 
$\e/\delta$ which appears in (\ref{piepsdelta}) is equal to 
$\lambda$. As $\e=\lambda\,\delta$, we can use the simplified notation
$u_i^\delta:=u_i^{\e,\delta}$ and $K_i^\delta:=K_i^{\e,\delta}$. Note 
that $(u_{0}^{\delta}, K_{0}^{\delta}):=(u_0, K_0)$ and 
for every $i\ge1$
$(u_i^{\delta}, K_i^{\delta})$ is a solution of 
the minimum problem
\begin{equation}\label{pidelta1}
\min_{(u,K)}\big\{ E(u,K)+
\lambda\,\|u- u_{i-1}^{\delta}\|^2
\big\}\,,
\end{equation}
where the minimum is taken over all pairs $(u,K)$ such that
$K\in\Kmfbul$, $K\supset K_{i-1}^{\delta}$, $u\in H^1(\Omk)$, 
$u=g_i^\delta$ on $\partial_D\Omk$. 

The term containing $\lambda$ is the main difference with respect to 
the discrete-time  
formulation of the model proposed by Francfort and Marigo in 
\cite{FraMar3}, which corresponds to the case $\lambda=0$. The fact 
that $\lambda$ is greater than $0$ penalizes the $L^2$-distance between 
$u_i^\delta$ 
and  $u_{i-1}^{\delta}$ and avoids some unnatural 
jumps which may occur in the continuous-time formulation for 
$\lambda=0$. In a sense, when $\lambda$ is large, local minimizers 
(close to $ u_{i-1}^{\delta}$) are preferred to global minimizers.

We introduce as before the piecewise constant interpolation 
$(u^\delta(t),K^\delta(t))$ defined by 
$u^{\delta}(t):=u_{i}^{\delta}$ and 
$K^{\delta}(t):=K_{i}^{\delta}$,
 for
$t_i^\delta\leq t<t_{i+1}^\delta$.

We prove (Lemma \ref{Helly2}) that there exists a left-continuous increasing 
function 
$K:[0,T]\to\Kmf$ such that for every $t\in[0,T]$, with 
$$
K(t)=K(t+):=\textstyle\bigcap_{s>t}K(s)\,,
$$
$K^\delta(t)$ converges to 
$K(t)$ as $\delta\to0$ along a suitable sequence  independent of $t$. 
In the rest of this discussion we always 
refer to this sequence when we write $\delta\to0$. 

For every $t\in[0,T]$ let $u(t)$ be a minimizer of 
\begin{equation}\label{et}
\int_{\Om\setminus K(t)}|\nabla u|^2\,dx 
\end{equation}
among all functions $u\in H^1(\Omk(t))$ with $u=g(t)$ on 
$\partial_D\Omk(t)$. We prove (Lemma~\ref{condc}) that 
\begin{equation}\label{b1}
E(u(t),K(t))\le E(u,K)+\lambda\,\|u-u(t)\|^2
\end{equation}    
for every $0<t\le T$,  for every $K\in\Kmf$ with $K\supset K(t)$, and 
for every $u\in H^1(\Omk)$ with $u=g(t)$ on $\partial_D\Omk$.
Moreover we prove (Lemma~\ref{ineq}) that
\begin{equation}\label{efg}
E(u(t),K(t))-E(u(s),K(s))\le 2\int_s^t\int_{\Om\setminus K(\tau)}\nabla 
u(\tau)\nabla \gdot(\tau)\,dx\,d\tau\,,
\end{equation}    
where $\gdot(t)$ is the time derivative of the function $g(t)$. 

This inequality shows that $t\mapsto E(u(t),K(t))$ is a function with
bounded variation and that 
$$
\frac{d}{dt}E(u(t),K(t))\le 2\int_{\Om\setminus K(t)}\nabla 
u(t)\nabla \gdot(t)\,dx
$$
for a.e.\  $t\in[0,T]$ (Remark~\ref{3.8}), and this leads to the existence of a function 
$\om(s,t)$, defined for $0\le s<t\le T$, with 
$$
\lim_{s\to t-}\frac{\om(s,t)}{t-s}=0\qquad\hbox{for every 
}t\in(0,T)\,,
$$
such that for a.e.\  $t\in[0,T]$ and every $s<t$ we have
\begin{equation}\label{om}
E(u(t),K(t))\le E(u,K(s))+\om(s,t)
 \end{equation}
for every $u\in H^1(\Omk(s))$ with $u=g(t)$ on $\partial_D\Omk(s)$
(Proposition~\ref{gtKs} and Remark~\ref{label}).

The minimality properties (\ref{b1}) and (\ref{om}) are used in 
Section~\ref{tips} to prove that the classical Griffith's criterion holds at 
the crack tips for a.e.\  $t\in[0,T]$, provided $K(t)$ satisfies suitable 
regularity conditions.

The fact that $\lambda >0$ in (\ref{pidelta1}) leads to an additional 
condition 
on the possible discontinuites of $(u(t),K(t))$. For every 
$t\in[0,T)$, for every $K\in\Kmf$,  and for every $u\in H^1(\Omk)$ 
with $u=g(t)$ on $\partial_D\Omk$, we determine a set ${\mathcal 
R}^t(u,K)$ (Definitions~\ref{r0} and~\ref{rt}), 
depending on the boundary condition $g(\cdot)$, such that for every
$t\in[0,T)$ we have $(u(t+),K(t+))\in {\mathcal 
R}^t(u(t),K(t))$, where $u(t+)$ is a minimizer of (\ref{et}) with $K(t)$ 
replaced by $K(t+)$ (Lemmas~\ref{condb} and~\ref{condd}).  
In Section~\ref{esempio} we show that, if $\Om$ and $g(t)$ are 
sufficiently regular and no initial crack is present, i.e., $K_0=\emptyset$, 
then, for $\lambda$ large enough, no crack will appear, i.e., 
$K(t)=\emptyset$ for every $t\in[0,T]$, and 
${\mathcal R}^t(u(t),\emptyset)=\{(u(t),\emptyset)\}$.
Note that the model by Francfort and Marigo \cite{FraMar3}, based on 
global minimization, gives, in general, a different result.

\end{section}

\begin{section}{NOTATION AND PRELIMINARIES}\label{notation}

Throughout the paper $\Om$ is a fixed bounded connected open subset
of $\R^2$ with Lipschitz
boundary.
Let $\K$ be the set of all compact subsets
of $\overline\Om$. Given an integer $m\ge 1$, let
$\Kmf$ be the set of all compact subsets $K$ of $\overline\Om$
with at most $m$ connected components and such that $\huno(K)<+\infty$.
We shall consider also the set $\Kmfbul$ of all $K\in\Kmf$ without 
isolated points.

We recall that the {\it Hausdorff distance\/} between
$K_{1},\, K_{2}\in\K$ is defined by
$$
d_{H}(K_{1},K_{2}):=
\max\big\{ \sup_{x\in K_1}{\rm dist}(x,K_2),
\sup_{y\in K_2}{\rm dist}(y,K_1)\big\}\,,
$$
with the conventions  ${\rm dist}(x,\emptyset)={\rm diam}(\Om)$ and
$\sup\emptyset=0$, so that $d_{H}(\emptyset, K)=0$ if $K=\emptyset$,
and $d_{H}(\emptyset, K)={\rm diam}(\Om)$ if $K\neq\emptyset$.
We say that $K_n\to K$ in the Hausdorff metric if
${d_H(K_n,K)\to0}$. The following compactness theorem is
well-known (see, e.g., 
\cite[Blaschke's Selection Theorem]{Rog}).
\begin{theorem}\label{compactness}
The metric space $(\K,d_H)$ is compact.
\end{theorem}
It is well-known that, in general, the Hausdorff measure is not
lower semicontinuous on $\K$
with respect to the convergence in the
Hausdorff metric. The following result, which is a consequence
of the Go\l \c ab theorem,
shows that on
the class $\Kmf$ the Hausdorff measure is lower semicontinuous. We
refer to \cite[Corollary~3.3]{DM-T} for a proof.
\begin{theorem}\label{Golab2}
Let $K_n$ be compact sets in $\overline\Om$
 with a uniformly bounded number of connected components. 
If $K_n$ converge to $K$ in the Hausdorff metric, then
$$
\huno(K\cap U)\le \liminf_{n\to\infty} \,\huno(K_n\cap U)
$$
for every open set $U\subset \R^2$. 
\end{theorem}

In the rest of the paper $\partial\Om$ is the union of
two (possibly empty) disjoint sets $\partial_D\Om$ and 
$\partial_N\Om$, with a finite
number of connected components,
on (part of) which we impose a nonhomogeneous 
Dirichlet boundary condition
and a homogeneous
Neumann
boundary condition, respectively.

Given a function $u\in H^1(\Omk)$ for some
$K\in\K$, we always
extend $u$ and $\nabla u$ to $\Om$ by setting $u=0$ and $\nabla u=0$
a.e.\ on $K$.
Note that, however, $\nabla u$ is the  distributional gradient  of $u$
only in $\Omk$, and, in general, $\nabla u$ does not coincide in $\Om$ 
with
the gradient of an extension of~$u$.

For every $K\in\K$ we consider the space
\begin{eqnarray*}
&& H^1_0(\Omk,\partial_D\Omk):=\{u\in H^1(\Omk):u=0
\hbox{ on }
\partial_D\Omk\}\,,
\end {eqnarray*}
where the equality on the boundary is intended in the sense of traces. 
The following definition reformulates in our particular case a general 
notion of convergence studied by Mosco in~\cite{Mos}.

\begin{definition}Let $K_n,K\in\K$. We say that the spaces
$H^1_0(\Omk_n,\partial_D\Omk_n)$ converge 
to the space $H^1_0(\Omk,\partial_D\Omk)$ in
the sense of Mosco if
the following properties hold:
\begin{itemize}
\item[\rm{($\mskip-.7\thinmuskip{\rm M}_1\mskip-1.3\thinmuskip$)}]
for every $u\in H^1_0(\Omk,\partial_D\Omk)$ there exists
$u_n\in H^1_0(\Omk_n,\partial_D\Omk_n)$ such that $u_n\to u$ strongly in
$L^2(\Om)$ and
$\nabla u_n\to \nabla u$ strongly in $L^2(\Om;\R^2)$;
\item[\rm{($\mskip-.7\thinmuskip{\rm M}_2\mskip-1.3\thinmuskip$)}] if
$u_n\in H^1_0(\Omk_{n},\partial_D\Omk_{n})$ for every $n$ and
$\liminf_{n}\|u_n\|_{H^1(\Om\setminus K_{n})}<+\infty$, then there exist a
subsequence $u_{n_k}$ and a function $u\in H^1_0(\Omk,\partial_D\Omk)$
such that
$u_{n_k}\wto u$ weakly in
$L^2(\Om)$ and $\nabla u_{n_k}\wto \nabla u$ weakly in
$L^2(\Om;\R^2)$.
\end{itemize}
\end{definition}

The following theorem shows the connection between Mosco convergence 
of the spaces  $H^1_0(\Omk_n,\partial_D\Omk_n)$ and convergence in 
the Hausdorff metric of the sets~$K_n$.

\begin{theorem}\label{Mosco} Let $K_n,K$ be compact sets in 
$\overline\Om$ with a uniformly bounded number of connected 
components, such that $K_n\to K$ in the Hausdorff metric and
${\rm meas}(K_n)\to{\rm meas}(K)$. Then
$H^1_0(\Omk_n,\partial_D\Omk_n)$
converges to $H^1_0(\Omk,\partial_D\Omk)$ in
the sense of Mosco.
\end{theorem}
\begin{proof} Under these hypotheses
Bucur and Varchon proved in \cite{BucVar} the Mosco convergence of
$H^1(\Omk_n)$ to
$H^1(\Omk)$. The extension  to the case when boundary
conditions are
imposed is due to Chambolle \cite{Ch} (see also \cite[Theorem~6.3]{DM-E-P}).
\end{proof}

Throughout the paper $\lambda$ is a fixed constant, with
$\lambda>0$. We use the notation
$(\cdot|\cdot)$  and $\|\cdot\|$
for the scalar product  and the norm in $L^2(\Om)$ or in
$L^2(\Om;\R^2)$,
according to the context. We have often to minimize energies of the 
form (\ref{energy}) among pairs $(u, K)$, where $K\in \Kmf$ and $u\in 
H^1(\Omk)$, with a prescribed boundary condition $u=\phi$ on 
$\partial_D\Omk$. We prefer to minimize first with respect to $u$ and 
then with respect to $K$. This leads to the following definitions.

\begin{definition}\label{def}
    If $K\in\Kmf$ for some $m\ge 1$, 
$\phi\in H^1(\Omk)\cap
L^\infty(\Om)$,  and
$w\in L^2(\Om)$,  we define
\begin{eqnarray}
&\displaystyle
\E(\phi,K):= \min_{v\in {\mathcal V}(\phi,K)}\{\|\nabla v\|^2+
\huno(K)\}\,,\label{egk}\\
&\displaystyle
\E_\lambda(\phi,K,w):= \min_{v\in {\mathcal V}(\phi,K)}
\{\|\nabla v\|^2+\huno(K)+\lambda\,\|v-w\|^2\}\,,\label{elambda}
\end{eqnarray}
where
\begin{equation}
{\mathcal V}(\phi,K):=\{v\in H^1(\Omk):v=\phi\quad\hbox{ on }
\partial_D\Omk\}\,.\label{vgk}
\end{equation}
\end{definition}

\begin{remark}\label{arm}
{\rm By minimality, a solution $u$ of (\ref{egk}) satisfies the 
inequality $\|\nabla u\|^2\le \|\nabla \phi\|^2$.
A truncation argument shows that there exists a minimizing
sequence $u_n$ of (\ref{egk}) such that $\|u_n\|_{\infty}\le
\|\phi\|_{\infty}$, where $\|\cdot\|_\infty$ denotes the norm in 
$L^\infty(\Om)$. By the direct method of the calculus of
variations we can then prove that there exists a solution $u$ of
(\ref{egk}) with
$\|u\|_{\infty}\le \|\phi\|_{\infty}$.  It is easy to
see that the solution is unique on the connected components of $\Omk$
whose boundaries meet $\partial_D\Omk$, while on the other connected
components it is given by an arbitrary constant. This shows
that two solutions have the same gradient.
If $u$ is a solution of the minimum problem
   (\ref{egk}), then $\E_\lambda(\phi,K,u)=\E(\phi,K)$. }
\end{remark}

   \begin{remark}\label{Linfty}{\rm By
   minimality, the solution
   $u$ of (\ref{elambda})  satisfies 
   $\|\nabla u\|^2+\lambda\,\|u-w\|^2\le \|\nabla
   \phi\|^2+\lambda\,\|\phi-w\|^2$.
   If  $w$
   belongs to $L^\infty(\Om)$, then an easy truncation argument
   shows that $u$
   belongs to $L^\infty(\Om)$  and
   $\|u\|_{\infty}\le \max
   \{\|\phi\|_{\infty}\,,\,\|w\|_{\infty}\}$. }
\end{remark}

\begin{remark}\label{costanti} {\rm If $w$ is constant on a connected
component $U$ of $\Omk$ whose boundary does not meet
$\partial_D\Omk$, then the minimizer $u$ of (\ref{elambda}) coincides
with $w$ on $U$. Therefore the value
$\E_\lambda(\phi,K,w)$ does not depend on the constant value of $w$ on
$U$.}
\end{remark}

\begin{remark}\label{equazione}{\rm
    A function $u$ is a minimizer of (\ref{egk}) if and only if
\begin{equation}\label{Pgk}
\left\{\begin{array}{ll} \Delta u=0 & \hbox{in }\Omk\,,\\
\salt
\frac{\partial u}{\partial\nu}=0 &
\hbox{on }
\partial_N\Om\cup K \,,\\
\salt
   u=\phi &\hbox{on } \partial_D\Om\setmeno K\,,
\end{array}\right.
\end{equation}
i.e., $u$ satisfies the following
conditions
\begin{equation}\label{P*gk}
\left\{\begin{array}{l}
u-\phi\in H^1_0(\Omk,\partial_D\Omk)\,,\\
\displaystyle\vphantom{\frac{\partial u}{}}
(\nabla u|\nabla  v)=0\qquad\forall
\,v\in H^1_{0}(\Omk,\partial_D\Omk)\,.
\end{array}\right.
\end{equation}
Similarly, $u$ is the minimizer of (\ref{elambda}) if and only if
$u$  is the solution of the problem
\begin{equation}\label{P}
\left\{\begin{array}{ll} \Delta u=\lambda(u-w) & \hbox{in }\Omk\,,\\
\salt
\frac{\partial u}{\partial\nu}=0 &
\hbox{on }
\partial_N\Om\cup K\,,\\
\salt
   u=\phi &\hbox{on } \partial_D\Om\setmeno K\,,
\end{array}\right.
\end{equation}
i.e., $u$ satisfies the following
conditions
\begin{equation}\label{P*}
\left\{\begin{array}{l}
u-\phi\in H^1_0(\Omk,\partial_D\Omk)\,,\\
\displaystyle\vphantom{\frac{\partial u}{}}
(\nabla u|\nabla  v)+\lambda(u-w|v)=0\qquad\forall
\,v\in H^1_{0}(\Omk,\partial_D\Omk)\,.
\end{array}\right.
\end{equation}
This implies that, if the minimizer $u$ of (\ref{elambda}) is equal 
to $w$, then $u$ is also a minimizer of~(\ref{egk}).
}
\end{remark}
We consider now the stability of the solutions to problems (\ref{egk}) 
and (\ref{elambda}) when $\phi$, $K$, and $w$ vary.

\begin{theorem}\label{C}
Let $m\ge1$, and let $\phi_n,\,\phi\in H^1(\Om)\cap L^\infty(\Om)$, 
$K_n,\,K\in\Kmf$, and $w_n,\,w\in 
L^2(\Om)$.  Let $u_n$ and $u$ be the solutions of the minimum problems 
(\ref{elambda}) which define $\E_\lambda(\phi_n,K_n,w_n)$ and 
$\E_\lambda(\phi,K,w)$, respectively. Assume that
$\phi_n\wto\phi$ weakly in $H^1(\Om)$, $K_n\to K$ in the Hausdorff 
metric, and $w_n\wto w$ weakly in $L^2(\Om)$. Then 
$u_n\wto u$ weakly in $L^2(\Om)$, $\nabla u_n\wto \nabla u$ 
weakly in $L^2(\Om;\R^2)$, and 
\begin{equation}\label{liminf}
\E_\lambda(\phi,K,w)\le\liminf_{h\to\infty}\E_\lambda(\phi_n,K_n,w_n)\,.
\end{equation}
 If $\phi_n$ and $w_n$ are uniformly bounded in 
$L^\infty(\Om)$ and $\phi_n\to\phi$ strongly in $H^1(\Om)$,  
then $u_n\to u$ strongly in $L^2(\Om)$ and
$\nabla u_n\to \nabla u$ 
strongly in $L^2(\Om;\R^2)$. 
\end{theorem}

\begin{proof} By Remark \ref{Linfty} the norms 
$\|u_n\|_{H^1({\Om\setminus K_n})}$ 
are uniformly bounded.
By Theorem \ref{Mosco} there exists $u^*\in H^1(\Omk)$, with $u^*=\phi$ 
on $\partial_D\Omk$, such that,  up to a subsequence, $u_n\wto u^*$ 
weakly in $L^2(\Om)$ and $\nabla u_n\wto \nabla u^*$   
weakly in $L^2(\Om;\R^2)$.
By (\ref{P*}) we have 
\begin{equation}\label{eh}
(\nabla u_n|\nabla  v_n)+\lambda(u_n-w_n|v_n)=0
\end{equation}
for every $v_n\in H^1_{0}(\Omk_n,\partial_D\Omk_n)$.
If $v\in H^1_{0}(\Omk,\partial_D\Omk)$, 
by Theorem~\ref{Mosco} there exist  
$v_n\in H^1_{0}(\Omk_n,\partial_D\Omk_n)$ such that $v_n\to v$ 
strongly in $L^2(\Om)$ and $\nabla v_n\to \nabla v$ 
strongly in $L^2(\Om;\R^2)$, and passing to the limit in (\ref{eh}) 
we obtain that $u^*$ is a solution of (\ref{P*}). By uniqueness, 
$u^*=u$, and the convergence holds for the whole sequence.
 Inequality (\ref{liminf}) follows now by lower semicontinuity 
 (Theorem~\ref{Golab2}).  

Assume that $\phi_n$ and $w_n$ are uniformly bounded in 
$L^\infty(\Om)$. Then the same is true for the solutions $u_n$
(Remark~\ref{Linfty}).
To prove the strong convergence in $L^2(\Om)$ of $u_n$, let
$U\subset\subset\Omk$ be an open set with boundary of class
$C^1$. As $K_n\to K$ in the Hausdorff
metric, we have $U\subset\subset\Omk_n$ for
$n$ large enough. Since
$u_n$ is  bounded in $H^1(U)$ uniformly with respect 
to~$n$ (Remark~\ref{Linfty}),
by the Rellich theorem
$u_n\to u$ strongly in $L^2(U)$.
As the functions $u_n$ are uniformly bounded in $L^\infty(\Om)$, 
the norms $\|u_n\|_{L^2(\Om\setminus U)}$ can be made arbitrarily small
by taking $U$ arbitrarily close to $\Omk$.
Therefore $u_n\to u$ strongly in $L^2(\Om)$.

If, in addition, $\phi_n\to\phi$ strongly in $H^1(\Om)$,
taking $v_n:=u_n-\phi_n$ as test function in (\ref{eh}) 
we can easily prove that
$\|\nabla u_n\|\to\|\nabla u\|$, which implies the strong convergence 
of the gradients.
\end{proof}
  
 The following corollary will be used in Section~\ref{proof}.
 
\begin{corollary}\label{CC}
Let $m\ge1$, and let $\phi_n,\,\phi\in H^1(\Om)\cap L^\infty(\Om)$, 
$K_n,\,K\in\Kmf$, and $w_n\in 
L^\infty(\Om)$. Let $u_n$ be the solutions of the minimum problems 
(\ref{elambda}) which define $\E_\lambda(\phi_n,K_n,w_n)$, and let
$u$ be a solution of the minimum problem (\ref{egk}) which 
defines $\E(\phi,K)$. 
Assume that $\phi_n$ and $w_n$ are uniformly bounded in 
$L^\infty(\Om)$, and that
$\phi_n\to\phi$ strongly in $H^1(\Om)$, $K_n\to K$ in the Hausdorff 
metric, and $u_n-w_n\to 0$ strongly in $L^2(\Om)$. Then
$\nabla u_n\to \nabla u$ 
strongly in $L^2(\Om;\R^2)$. Moreover there exist 
a subsequence $u_{n_k}$ of $u_n$ and a solution $u^*$ 
of the minimum problem (\ref{egk}) which 
defines $\E(\phi,K)$ (possibly different from~$u$) such that
$u_{n_k}\to u^*$ strongly in $L^2(\Om)$.
\end{corollary}

\begin{proof}
As $w_n$ is bounded in $L^\infty(\Om)$, there exists a subsequence 
$w_{n_k}$ which converges weakly in $L^2(\Om)$ to a function $w$.
Let $u^*$ be the solution of the minimum problem
(\ref{elambda}) which defines $\E_\lambda(\phi,K,w)$.
By Theorem~\ref{C} we have $u_{n_k}\to u^*$
strongly in $L^2(\Om)$ and $\nabla u_{n_k}\to \nabla u^*$ strongly 
in $L^2(\Om;\R^2)$. As 
$u_n-w_n\to 0$ strongly in $L^2(\Om)$, the functions $u^*$ and $w$
are equal. By 
Remark~\ref{equazione} this implies that $u^*$
is a solution of the minimum problem 
(\ref{egk}) which defines $\E(\phi,K)$, and
by the uniqueness of the gradients
we have $\nabla u=\nabla u^*$ a.e.\ in $\Om$.
Since we can repeat this argument for 
an arbitrary subsequence, we conclude that the
whole sequence
$\nabla u_n$ converges to $\nabla u$ strongly in
$L^2(\Om;\R^2)$.
\end{proof}
   
\end{section}

\begin{section}{IRREVERSIBLE QUASI-STATIC EVOLUTION}
\label{main}

In this section we define a continuous-time evolution of a cracked body
by investigating the properties of the
limits of the discrete-time evolution described in the introduction.

Let us fix the boundary displacement
$g\in AC([0,T];H^1(\Om))\cap L^\infty([0,T];L^\infty(\Om))$
and an integer $m\ge 1$.
Given an initial crack $K_0\in\Kmfbul$, we shall construct an
   increasing function $ K\colon [0,T]\to\Kmf$ satisfying suitable
   minimality conditions.
We define
\begin{eqnarray}
\textstyle K(t-):={\rm cl}\big(\bigcup_{s<t} K(s)\big)\
&&\hbox{ for }0<t\leq T\,,
\label{k_*t}\\
\textstyle  K(t+):=\bigcap_{s>t}K(s)\qquad &&\hbox{ for }0
\leq t<T\,,\label{k^*t}
\end{eqnarray}
where ${\rm cl}$ denotes the closure. We say that $t\mapsto K(t)$ is 
left-continuous if $K(t-)=K(t)$ for every $t\in (0,T]$.
It is easy to see that
\begin{eqnarray}
K(t-)\subset  K(t)\subset  K(t+)\,\,\quad
&&\hbox{ for }0< t<T\,,\label{7.9}\\
\textstyle K(t-)={\rm cl}\big(\bigcup_{s<t} K(s-)\big)\
&&\hbox{ for }0<t\leq T\,,\label{k-}
\\
\textstyle  K(t+)=\bigcap_{s>t}K(s+)\qquad &&\hbox{ for }0
\leq t<T\,.\label{k+}
\end{eqnarray}
Let $\Theta$ be the set of points $t\in (0,T)$ such that $K(t-)=K(t+)$.
By \cite[Theorem~6.1]{DM-T} the set $[0,T]\setmeno\Theta$ is at most 
countable.

For every $t\in[0,T]$ let $u(t)$ (resp.\  $u(t-)$, $u(t+)$) be a solution of
the minimum problem (\ref{egk}) corresponding to $\phi =g(t)$ 
and $K=K(t)$ (resp.\  $K=K(t-)$, $K=K(t+)$). By Remark~\ref{arm} 
$\|\nabla u(t)\|$ is bounded uniformly with respect to~$t$.
By \cite[Theorem~5.1]{DM-T} and by (\ref{k_*t})
and (\ref{k^*t}),
\begin{eqnarray}
\hbox{ for }0<t\le T&& \nabla u(s)\to  \nabla u(t-)  \qquad\hbox{strongly in
}L^2(\Om;\R^2)\hbox{ as }s\to {t-}\,,\label{ut-}\\
   \hbox{ for }0\le t<T && \nabla u(s)\to  \nabla u(t+)
\qquad\hbox{strongly in
}L^2(\Om;\R^2)\hbox{ as }s\to {t+}\,.\label{ut+}
\end{eqnarray}
This implies in particular that $t\mapsto \nabla u(t)$ is continuous 
from $[0,T]$ into $L^2(\Om;\R^2)$ at every point $t\in \Theta$. 
Therefore the first estimate of Remark~\ref{arm} implies that
\begin{equation}\label{measurab}
t\mapsto \nabla u(t)
\quad\hbox{belongs to}\quad
L^\infty([0,T];L^2(\Om;\R^2))\,.
\end{equation}

Although the boundary displacement $g(t)$ is continuous with respect to
$t$, the continuous-time evolution that we shall obtain as limit of the
discrete-time evolutions may exhibit some jump discontinuities of the pair
$(u(t),K(t))$. 
Given a time step $\delta>0$, the approximation procedure considered 
in the introduction uses 
sequences $(v_i^\delta,H_i^\delta)$ with the property that, for every 
$i\ge1$, $H_i^\delta$ is a solution of the minimum problem
\begin{equation}\label{pidelta3}
\min_K\big\{ \E_\lambda(g_i^\delta,K,v_{i-1}^\delta) : K\in\Kmfbul,
\ K\supset H_{i-1}^\delta\big\}\,,
\end{equation}
and $v_i^\delta$ is the solution of the minimum problem
(\ref{elambda})  defining
$\E_\lambda(g_i^\delta,H_i^\delta,v_{i-1}^\delta)$. We recall that 
$g_i^\delta:=g(t_i^\delta)$, with $t_i^\delta:=i\delta$. The existence of 
a solution to (\ref{pidelta3}) is proved in Lemma~\ref{3.5}.

Let us consider first the discontinuities that may occur at the
initial time  $t=0$.
\begin{definition}\label{r0}
Given a pair $(u,K)$ with $K\in\Kmfbul$, $u\in H^1(\Omk)$, 
and $u=g(0)$ on $\partial_D\Omk$,
we  define ${\mathcal R}^0(u,K)$
as the set of  all pairs $(v,H)$ such that
\begin{itemize}
    \item[(a)] $H\in\Kmf$,
    $v\in H^1(\Om\setmeno H)$, $v=g(0)$ on 
    $\partial_D\Om\setmeno H$,
   \item[(b)]  there exist a sequence $\delta_n\to 0+$, a 
   sequence of integers 
   $l_n\to\infty$, with $l_n\delta_n\to0$, and a sequence
   $(v_i^\delta,H_i^\delta)$ satisfying (\ref{pidelta3}), such that
  \begin{itemize}
      \item[$({\rm b}_1)$]
      $\vphantom{ H_{l_n}^{\delta_n}}v_0^{\delta_n}=u$ 
      and $H_0^{\delta_n}=K$ for every $n$,
      \item[$({\rm b}_2)$] $\vphantom{ H_{l_n}^{\delta_n}}H_{l_n}^{\delta_n}\to 
      H$ in the Hausdorff metric and $\vphantom{ H_{l_n}^{\delta_n}}\nabla 
      v_{l_n}^{\delta_n}\to\nabla v$ strongly in $L^2(\Om;\R^2)$.
      \end{itemize}
  \end{itemize}
  \end{definition}
We will prove that the continuous-time evolution  satisfies
$$
( u(0+),K(0+))\in {\mathcal R}^0(u(0),K(0))\,.
$$
This shows in particular that $(u(0+),K(0+))=(u(0),K(0))$ when
${\mathcal R}^0(u(0),K(0))$
contains only $(u(0),K(0))$ (see Section~\ref{esempio}).

The definition of ${\mathcal R}^t(u,K)$ at time $t>0$ is more complex, since
the approximation procedure described in the introduction forces us to
replace $({\rm b}_1)$ by a  more technical
condition.
\begin{definition}\label{rt}
Given $t\in(0,T)$ and a pair $(u,K)$ with $K\in\Kmfbul$, $u\in H^1(\Omk)$, 
and $u=g(t)$ on $\partial_D\Omk$, 
we define ${\mathcal R}^t(u,K)$ as  the set of 
all pairs $(v,H)$ such that
\begin{itemize}
    \item[(a)] $H\in\Kmf$,
    $v\in H^1(\Om\setmeno H)$, $v=g(t)$ on 
    $\partial_D\Om\setmeno H$,
   \item[(b)]  there exist a sequence $\delta_n\to 0+$, three 
   sequences of integers $h_n$, $k_n$, $l_n$ converging to~$\infty$, 
   with $h_n\delta_n\to t-$, $k_n\delta_n\to 
   t-$, $l_n\delta_n\to t+$, $k_n-h_n\to\infty$,
   $l_n-k_n\to\infty$, and a sequence
   $(v_i^\delta,H_i^\delta)$ satisfying (\ref{pidelta3}), such 
   that 
  \begin{itemize}
       \item[$({\rm b}_1)$] for every sequence 
       $\sigma_n$ with $h_n\le \sigma_n\le 
       k_n$ we have 
    $ H_{\sigma_n}^{\delta_n}\to K$ in the Hausdorff 
     metric and $\nabla  v_{\sigma_n}^{\delta_n}\to \nabla u$ strongly in 
$L^2(\Om;\R^2)$,
    \item[$({\rm b}_2)$] $\vphantom{ H_{l_n}^{\delta_n}}H_{l_n}^{\delta_n}\to 
      H$ in the Hausdorff metric and $\vphantom{ H_{l_n}^{\delta_n}}\nabla 
      v_{l_n}^{\delta_n}\to\nabla v$ strongly in $L^2(\Om;\R^2)$.
      \end{itemize}
  \end{itemize}
\end{definition}
  We will prove that the continuous-time evolution  satisfies
$$
( u(t+),K(t+))\in {\mathcal R}^t( u(t-),K(t-))\,.
$$
This gives a restriction on the possible jumps
and shows in particular that
$( u(t),K(t))$ is continuous at time $t$ whenever
${\mathcal R}^t( u(t-),K(t-))$  contains only $( u(t-),K(t-))$ (see
Section~\ref{esempio}).

We are now in a position to state the main result of the paper, which 
provides a continuous-time variational model for the 
quasi-static growth of brittle 
fractures.

\begin{theorem}\label{kt}
Let $T>0$, $\lambda>0$, $m\ge1$, let $g\in AC([0,T];H^1(\Om))\cap
L^\infty([0,T];L^\infty(\Om))$, let $\gdot\in L^1([0,T];H^1(\Om))$ be 
its time derivative, and let 
$K_{0}\in\Kmfbul$.
Then there exists a function
$K\colon[0,T]\to\Kmf$  such that, if $u(t)$ is a solution of the
minimum problem (\ref{egk})  corresponding to $\phi=g(t)$ and $K=K(t)$,
the following conditions are satisfied:
\smallskip
\begin{itemize}
\item[\rm{(a)}]
\hfil $\displaystyle \vphantom{\frac{d}{ds}}
K(0)= K_0  \quad\hbox{ and }\quad K(s)\subset K(t)\quad\hbox{for }0\le 
s\le t\le T\,,$\hfil
\item[\rm{(b)}] \hfil $\displaystyle \vphantom{\frac{d}{ds}}
\hbox{for  } 0< t \le
T\quad\E_\lambda(g(t),K(t),u(t))\leq
\E_\lambda(g(t),K,u(t))
\quad\forall \, K\in\Kmf\,,\,\  K\supset K(t)\,,$\hfil
\item[\rm{(c)}] \hfil$\displaystyle \vphantom{\frac{d}{ds}}\hbox{ for } 
0\le s\le t\le T\quad\E(g(t),K(t))-\E(g(s),K(s))\le
2\int_s^t(\nabla u(\tau)|\nabla \gdot(\tau))\,d\tau\,,$\hfil
\item[\rm{(d)}] \hfil $\displaystyle \vphantom{\frac{d}{ds}}
(u(0+),K(0+))\in {\mathcal R}^0(u(0),K(0))\,,$\hfil
\item[\rm{(e)}] \hfil $\displaystyle \vphantom{\frac{d}{ds}}
\hbox{for  } 0< t < T\quad (u(t+),K(t+))\in {\mathcal
R}^t(u(t-),K(t-))\,,$\hfil
\end{itemize}
where $K(t-)$ and $K(t+)$  are defined by (\ref{k_*t}) and 
(\ref{k^*t}), while $u(t-)$ and $u(t+)$ are solutions of the minimum 
problems (\ref{egk}), with $\phi=g(t)$, corresponding to $K=K(t-)$ 
and $K=K(t+)$.
\end{theorem}

\begin{remark}\label{3.7}{\rm Since $t\mapsto \nabla u(t)$ belongs to 
$L^\infty([0,T];L^2(\Om;\R^2))$ (see (\ref{measurab})) and 
$t\mapsto \nabla \gdot(t)$ belongs to 
$L^1([0,T];L^2(\Om;\R^2))$, the function $t\mapsto (\nabla u(t)|\nabla 
\gdot(t))$ is integrable on $[0,T]$.
}
\end{remark}

\begin{remark}\label{3.8} {\rm  By Remark~\ref{3.7},
condition (c) of Theorem~\ref{kt} implies that
\begin{itemize}    
\item[\rm{(f)}]   the
function $\displaystyle \vphantom{\frac{d}{}}
t\mapsto \E(g(t),K(t))$
   has bounded variation on $[0,T]$,
   and its positive variation is absolutely continuous on
$\displaystyle \vphantom{\frac{}{ds}}
[0,T]$;
\item[\rm{(g)}] \hfil$\displaystyle\frac{d}{dt}\E(g(t),K(t))\le
2(\nabla u(t)|\nabla \gdot(t))\quad \hbox{ for a.e.\ } t\in[0,T]$.\hfil
\end{itemize}
Conversely, (c) follows from (f) and (g). 
}
\end{remark}

\begin{proposition} \label{gtKs}
 Under the assumptions of Theorem~\ref{kt}, if 
$K:[0,T]\to \Kmf$ satisfies (a) and (c), then 
\smallskip
\begin{itemize}    
\item[\rm{(h)}]\hfil $\displaystyle \vphantom{\frac{d}{}}  
\lim_{s\to t}\frac{\E(g(t),K(s))-\E(g(t),K(t))}{s-t}\le 0\quad
\hbox{for a.e.\ }t\in[0,T]$.\hfil
\end{itemize}
\smallskip\noindent
Conversely, if $t\mapsto K(t)$ satisfies (a) of Theorem~\ref{kt}, (f) of 
Remark~\ref{3.8}, and
\smallskip
\begin{itemize}    
\item[\rm{($ \mskip-.7\thinmuskip{\rm h}'\mskip-.7\thinmuskip$)}]
\hfil $\displaystyle \vphantom{\frac{d}{}}  
\limsup_{s\to t-}\frac{\E(g(t),K(s))-\E(g(t),K(t))}{s-t}\le 0\quad
\hbox{for a.e.\ }t\in[0,T]$.\hfil
\end{itemize}
\smallskip\noindent
then $t\mapsto K(t)$ satisfies also (g) of 
Remark~\ref{3.8}; therefore it satisfies condition (c) of Theorem~\ref{kt}.
\end{proposition}
\begin{proof}
 We notice that $\E(g(t),K(t))-\E(g(s),K(s))$ can be written as
\begin{equation}\label{e1}
[\E(g(t),K(t))-\E(g(t),K(s))]+
[\E(g(t),K(s))-\E(g(s),K(s))]\,.
\end{equation}
Let $u(t,s)$ be a solution of the minimum problem (\ref{egk})
corresponding to $\phi=g(t)$ and $K=K(s)$. Then taking
$u(t,s)-u(s)-g(t)+g(s)$ as test function in the equations
satisfied by $u(t,s)$ and $u(s)$ we obtain that
\begin{equation}\label{e2}
\E(g(t),K(s))-\E(g(s),K(s))=
(\nabla u(t,s)+\nabla u(s)|\nabla g(t)-\nabla g(s))\,.
\end{equation}
Let $\Theta$ be the set of points $t\in(0,T)$ such that $K(t+)=K(t-)$. 
By \cite[Proposition 6.1]{DM-T} we have that $[0,T]\setmeno\Theta$ 
is at most countable. 
Assume that $t\in\Theta$. As $K(s)\to K(t)$ in the Hausdorff metric for
$s\to t$,
by \cite[Theorem~5.1]{DM-T} both $\nabla u(t,s)$ and 
$\nabla u(s)$ converge to
$\nabla u(t)$
strongly in $L^2(\Om;\R^2)$ as $s\to t$.
We now divide (\ref{e1}) and (\ref{e2}) by $t-s$ and pass
to the limit as $s\to t-$. If (c) is satisfied,
from condition (g) of Remark~\ref{3.8} we get 
(h) for all $t\in\Theta$ such that
$\frac{d}{dt}\E(g(t),K(t))$ and $\nabla \gdot(t)$ exist.
Conversely, if  (f) and
($ \mskip-.7\thinmuskip{\rm h}'\mskip-.7\thinmuskip$)
 are satisfied, then (g) holds
for all $t\in\Theta$ such that
$\frac{d}{dt}\E(g(t),K(t))$ and $\nabla \gdot(t)$ exist.
\end{proof}   
    
\begin{remark}\label{label}{\rm 
Condition ($ \mskip-.7\thinmuskip{\rm h}'\mskip-.7\thinmuskip$)
of Proposition \ref{gtKs} 
is equivalent to the existence of a
function $\om(s,t)$, defined for $0\le s<t\le T$, with 
$$
\lim_{s\to t-}\frac{\om(s,t)}{t-s}=0\qquad\hbox{for every 
}t\in(0,T)\,,
$$
such that for a.e.\  $t\in[0,T]$ and every $s<t$ the energy $E(u,K)$ 
defined in (\ref{energy}) satisfies
\begin{equation}\nonumber
E(u(t),K(t))\le E(u,K(s))+\om(s,t)
\end{equation}
for every $u\in H^1(\Omk(s))$ with $u=g(t)$ on $\partial_D\Omk(s)$.
}
\end{remark}
\end{section}

\begin{section}{PROOF OF THE MAIN RESULT}\label{proof}
In this section we prove
Theorem~\ref{kt} by a time discretization process.
Let us fix a solution
$u_0\in H^1({\Omk_0})\cap L^\infty(\Om)$ of the minimum
problem (\ref{egk})
corresponding to $\phi=g(0)$ and to $K=K_0$. By Remark~\ref{arm} we may 
assume that $\|u_0\|_\infty\le \|g(0)\|_\infty$.
Given $\delta>0$, we define $(u_i^\delta, K_i^\delta)$ inductively as follows:
   $u_{0}^\delta:=u_0$ and
$K_{0}^\delta:=K_0$;
for $i\geq1$ we define $K_i^\delta$ as a solution of the minimum problem
\begin{equation}\label{pidelta}
\min_K\big\{ \E_\lambda(g_i^\delta,K,u_{i-1}^\delta) : K\in\Kmfbul,
\ K\supset K_{i-1}^\delta\big\}\,,
\end{equation}
and $u_i^\delta$ as the solution of the minimum problem
(\ref{elambda})  defining
$\E_\lambda(g_i^\delta,K_i^\delta,u_{i-1}^\delta)$.

\begin{lemma}\label{3.5}
There exists a solution $K_i^\delta$ of the minimum problem (\ref{pidelta}). 
Moreover 
\begin{equation}\label{Kbul}
    \E_\lambda(g_i^\delta,K_i^\delta,u_{i-1}^\delta)\le 
    \E_\lambda(g_i^\delta,K,u_{i-1}^\delta)
\end{equation}
for every $K\in\Kmf$ with $K\supset K_{i-1}^\delta$.
\end{lemma}
\begin{proof}
By hypothesis $K_{0}^\delta:=K_0\in\Kmfbul$. Assume by induction that
$K_{i-1}^\delta\in \Kmfbul$ and $u_{i-1}^\delta\in
H^1(\Omk_{i-1}^\delta)\cap L^\infty(\Om)$.
Consider a minimizing sequence $K_n$ of problem~(\ref{pidelta}).
We may assume  that $\huno(K_n)$ is uniformly 
bounded.
By compactness (Theorem~\ref{compactness}), passing to a subsequence
we may assume that  $K_n$ converges in the Hausdorff metric  to
some  compact set $K^*$ containing $K_{i-1}^\delta$.  By 
Theorem~\ref{Golab2} we have  $K^*\in\Kmf$. Let $K_i^\delta$ be the set 
of nonisolated points of $K^*$.  Then $K_i^\delta\in\Kmfbul$ and 
$K^*\setmeno K_i^\delta$ has a finite number of points.
Since $K^*\supset K_{i-1}^\delta$ and $K_{i-1}^\delta$ has no isolated points, 
we have
$K_i^\delta\supset K_{i-1}^\delta$.
By Theorem~\ref{C}
we conclude that
$\E_\lambda(g_i^\delta,K_i^\delta,u_{i-1}^\delta)=
\E_\lambda(g_i^\delta,K^*,u_{i-1}^\delta)\leq
\liminf_{n}\E_\lambda(g_i^\delta,K_n,u_{i-1}^\delta)$.
Since
$K_n$ is a minimizing sequence, this proves that $K_i^\delta$  is a
solution of the minimum problem (\ref{pidelta}). 

To prove (\ref{Kbul}) it is enough to observe that if  $K\in\Kmf$ and 
$K\supset K_{i-1}^\delta$, the set $K'$ of nonisolated points of $K$ 
belongs to $\Kmfbul$ and contains $K_{i-1}^\delta$ (since this set has 
no isolated points). As $K\setmeno K'$ has a finite number of points, 
from (\ref{pidelta})  
we obtain 
$\E_\lambda(g_i^\delta,K_i^\delta,u_{i-1}^\delta)\le
\E_\lambda(g_i^\delta,K',u_{i-1}^\delta)=
\E_\lambda(g_i^\delta,K,u_{i-1}^\delta)$.
\end{proof}

\begin{remark}\label{linfty} {\rm
   If $M$ is a constant such that
   $\|g(t)\|_{\infty}\le M$ for every $t\in[0,T]$, then
   $\|u_0\|_{\infty}\le M$ and
   Remark~\ref{Linfty}, applied inductively, gives
   $\|u_i^\delta\|_{\infty}\le M$ for every $\delta>0$ and every
   $i\ge 0$ with $t_i^\delta\le T$.
   By the minimality of $u_i^\delta$ we have
$$
   \|\nabla u_i^\delta\|^2\le 
   \|\nabla g_i^\delta\|^2 + \lambda\,\|g_i^\delta-u_{i-1}^\delta\|^2\,,
   $$
   which shows that $\nabla u_i^\delta$ is bounded in $L^2(\Om,\R^2)$
   uniformly with respect to $\delta$ and $i$.}
\end{remark}
We define now the step functions $g^\delta(t)$, $K^\delta(t)$, and 
$u^\delta(t)$
on $[0,T]$
by setting $g^\delta(t):=g_{i}^\delta$, $K^\delta(t):=K_{i}^\delta$,
and $u^\delta(t):=u_{i}^\delta$ for
$t_i^\delta\leq t<t_{i+1}^\delta$.

\begin{lemma}\label{discr}
There exists a positive function $\rho(\delta)$, converging to zero
as $\delta\to0$, such that
\begin{eqnarray}\label{2discr}
&
\|\nabla u_j^\delta\|^2+\huno(K_j^\delta)+
\lambda\displaystyle\sum_{h=i+1}^j\|u_h^\delta-u_{h-1}^\delta\|^2
\leq\\\nonumber
&\leq\|\nabla u_i^\delta\|^2+\huno(K_i^\delta)+
2\displaystyle\int_{t_i^{\delta}}^{t_j^{\delta}}(\nabla u^\delta(\tau)|\nabla
\gdot(\tau))\,d\tau+\rho(\delta)
\end{eqnarray}
for $0\leq i<j$ with $t_j^\delta\le T$.
\end{lemma}
\begin{proof}
Let us fix an integer $r$ with $i\leq r<j$. {}From the absolute
continuity of $t\mapsto g(t)$ we have
\begin{equation}\label{gr}
g_{r+1}^\delta-g_r^\delta=\int_{t_r^{\delta}}^{t_{r+1}^{\delta}}
\gdot(\tau)\,d\tau\,,
\end{equation}
where the integral is a Bochner integral for functions with values in
$H^1(\Om)$. This implies that
\begin{equation}\label{nabla}
\nabla g_{r+1}^\delta-\nabla g_r^\delta=
\int_{t_r^{\delta}}^{t_{r+1}^{\delta}}
\nabla \gdot(\tau)\,d\tau\,,
\end{equation}
where the integral is now a Bochner integral for functions with values in
$L^2(\Om;\R^2)$.

As $u_r^{\delta}+g_{r+1}^\delta-g_r^\delta\in H^1(\Omk_{r}^\delta)$
and $u_r^{\delta}+g_{r+1}^\delta-g_r^\delta=g_{r+1}^\delta$ on
$\partial_D\Omk_r^\delta$, we have
\begin{equation}\label{a}
\E_\lambda(g_{r+1}^\delta,K_r^\delta,u_r^{\delta} )\le
\|\nabla u_r^\delta+\nabla g_{r+1}^\delta-\nabla g_r^\delta\|^2+
\huno(K_r^\delta)+\lambda\,\|g_{r+1}^\delta-g_r^\delta\|^2\,.
\end{equation}
By the minimality of $u_{r+1}^\delta$ and of $K_{r+1}^\delta$ (see
(\ref{pidelta}))
we have
\begin{equation}\label{b}
\|\nabla
u_{r+1}^\delta\|^2+\huno(K_{r+1}^\delta)+\lambda\,\|u_{r+1}^\delta-u_r^\delta\|^2=
\E_\lambda(g_{r+1}^\delta,K_{r+1}^\delta,u_r^{\delta})\le
\E_\lambda(g_{r+1}^\delta,K_r^\delta,u_r^{\delta})\,.
\end{equation}
{}From (\ref{gr}), (\ref{nabla}), (\ref{a}), and (\ref{b})  we obtain
\begin{eqnarray*}
& \displaystyle \vphantom{\int_{t_r^{\delta}}^{t_{r+1}^{\delta}}
(\nabla u^\delta(\tau)|\nabla\gdot(\tau))\,d\tau}
\|\nabla u_{r+1}^\delta\|^2+\huno(K_{r+1}^\delta)+
\lambda\,\|u_{r+1}^\delta-u_r^\delta\|^2\leq\\
& \displaystyle \vphantom{\int_{t_r^{\delta}}^{t_{r+1}^{\delta}}
(\nabla u^\delta(\tau)|\nabla\gdot(\tau))\,d\tau}
\leq
\|\nabla u_r^\delta+\nabla g_{r+1}^\delta-\nabla g_r^\delta\|^2+
\huno(K_r^\delta)+\lambda\,\|g_{r+1}^\delta-g_r^\delta\|^2\leq\\
& \leq  \|\nabla u_{r}^\delta\|^2+\huno(K_r^\delta)+2
\displaystyle\int_{t_r^{\delta}}^{t_{r+1}^{\delta}}
(\nabla u_{r}^\delta|\nabla
\gdot(\tau))\,d\tau+\\
&{} +\Big(\displaystyle\int_{t_r^{\delta}}^{t_{r+1}^{\delta}}\|\nabla
\gdot(\tau)\|\,d\tau\Big)^2+\lambda
\Big(\displaystyle\int_{t_r^{\delta}}^{t_{r+1
}^{\delta}}\|
\gdot(\tau)\|\,d\tau\Big)^2\leq\\
&  \leq  \|\nabla u_{r}^\delta\|^2+\huno(K_r^\delta)+2
\displaystyle\int_{t_r^{\delta}}^{t_{r+1}^{\delta}}
(\nabla u^\delta(\tau)|\nabla\gdot(\tau))\,d\tau+\\
&{} +
\sigma(\delta)\Big(\displaystyle
\int_{t_r^{\delta}}^{t_{r+1}^{\delta}}\|\nabla
\gdot(\tau)\|\,d\tau+\lambda\displaystyle
\int_{t_r^{\delta}}^{t_{r+1}^{\delta}}
\|\gdot(\tau)\|\,d\tau\Big)\,,
\end{eqnarray*}
where
$$
\sigma(\delta):=\max_{r}
\int_{t_r^{\delta}}^{t_{r+1}^{\delta}}(\|\nabla
\gdot(\tau)\|+\|\gdot(\tau)\|)\,d\tau \ \longrightarrow \ 0
$$
by the absolute continuity of the integral.
Iterating now this inequality for $i\leq r<j$ we get
(\ref{2discr}) with $\rho(\delta):=\sigma(\delta)\int_0^T(\|\nabla
\gdot(\tau)\|+\lambda\,\|\gdot(\tau)\|)\,d\tau$.
\end{proof}

\begin{lemma}\label{estim}
There exists a constant $M$, depending only on $g$, $K_0$, $u_0$, and
$\lambda$,
such that
\begin{equation}\label{stima}
\|\nabla u_i^\delta\|^2\leq M\,,\quad
\sum_{0<t_h^\delta\le T}
\|u_h^\delta-u_{h-1}^\delta\|^2\leq M\,,
\quad\hbox{and}\quad
\huno(K_i^{\delta})\leq M
\end{equation}
for every $\delta>0$ and for every  $i\ge0$ with $t_i^\delta\le T$.
\end{lemma}
\begin{proof}
{}From the previous lemma we get
\begin{eqnarray}\nonumber
&
\|\nabla u_i^\delta\|^2+\huno(K_i^\delta)+
\lambda\displaystyle\sum_{h=1}^i\|u_h^\delta-u_{h-1}^\delta\|^2
\leq\\
\label{idiscr}
&\leq\|\nabla u_0^\delta\|^2+\huno(K_0^\delta)+
2\displaystyle\int_{t_0^{\delta}}^{t_i^{\delta}}
(\nabla u^\delta(\tau)|\nabla
\gdot(\tau))\,d\tau+\rho(\delta)=\\
\nonumber
&
=\|\nabla u_0\|^2+\huno(K_0)+
2\displaystyle\int_{0}^{t_i^{\delta}}(\nabla u^\delta(\tau)|\nabla
\gdot(\tau))\,d\tau+\rho(\delta)\,.
\end{eqnarray}
The first
inequality in (\ref{stima}) is proved in Remark~\ref{linfty}. The other
inequalities follow now from (\ref{idiscr}).
\end{proof}

\begin{lemma}\label{Helly2} 
There exists an increasing left-continuous
function $K\colon[0,T]\to\Kmf$
such that,
for every $t\in(0,T)$ with $K(t)=K(t+)$, 
$K^\delta(t)$ converges to $K(t)$ in the
Hausdorff metric as $\delta\to0$ along a suitable sequence
independent of~$t$.
\end{lemma}
\begin{proof}
By \cite[Theorem 6.3]{DM-T} there exists an increasing function
$\hat K\colon[0,T]\to\K$ such that,
for every $t\in[0,T]$, $K^\delta(t)$ converges to $\hat K(t)$ in the
Hausdorff metric as $\delta\to0$ along a suitable sequence
independent of~$t$. By
Lemma~\ref{estim} we have $\huno(K^\delta(t))\le M$ for every
$t\in[0,T]$, $\delta>0$, and by Theorem~\ref{Golab2} this
implies $\hat K(t)\in\Kmf$ for every $t\in[0,T]$.
Let $K\colon[0,T]\to\Kmf$ be
the left-continuous  regularization of $\hat K(t)$, 
defined by $K(0)=\hat K(0)$ and
$K(t):=\hat K(t-)$ for every $t\in(0, T]$. Then $K(t)$ is 
left-continuous by (\ref{k-}), and by (\ref{7.9})
$K^\delta(t)$ converges to $K(t)$ for every 
$t\in(0,T)$ with $K(t)=K(t+)$.
\end{proof}

In the rest of this section, when
we write $\delta\to0$ we always refer to the sequence given
by Lemma~\ref{Helly2}.
Let $\Theta$ be the set of points $t\in(0,T)$
such that  $K(t)=K(t+)$.
Then ${[0,T]\setmeno\Theta}$ is at most countable (see 
\cite[Proposition 6.1]{DM-T}), and
$K(t_n)\to K(t)$ in the Hausdorff metric for every $t\in\Theta$ and
every sequence $t_n$ in $[0,T]$ converging to~$t$.

\begin{lemma}\label{l2}
For every $t\in(0,T]$ there exist  two sequences of integers
$h_\delta$ and $k_\delta$
such that
$k_\delta-h_\delta\to\infty$,
$h_\delta\delta\to t-$, $k_\delta\delta\to t-$,
$K^\delta(h_\delta\delta)$ and  $K^\delta(k_\delta\delta)$
converge to
$K(t)$ in the Hausdorff metric,
and
\begin{equation}
\sum_{h=h_\delta}^{k_\delta}
\|u_h^\delta-u_{h-1}^\delta\|^2\to0\label{somma}
\end{equation}
as $\delta\to0$.
In particular, setting $t_\delta=h_\delta\delta$, we have that
$u^\delta(t_\delta)-u^\delta(t_\delta-\delta)\to0$ strongly in~$L^2(\Om)$.
\end{lemma}
\begin{proof} Let $\tau_k\to t-$ be such that
$\tau_k\in\Theta$. Then
$K^\delta(\tau_k)\to K(\tau_k)$ in the Hausdorff metric
as $\delta\to 0$
by Lemma~\ref{Helly2}. We choose a
strictly
decreasing  sequence $\delta_k\searrow0$ such that
for every
$\delta\le\delta_k$
$$
d_H(K^\delta(\tau_k),K(\tau_k))<\frac{1}{k}\,.
$$
For $\delta_{k+1}\le\delta<\delta_k$ let $s_\delta=\tau_k$.
Then $s_\delta\to t-$ and
$$
d_H(K^\delta(s_\delta),K(t))\leq
\frac{1}{k}+d_H(K(s_\delta),K(t))\qquad\hbox{for
}\quad\delta_{k+1}\le\delta<\delta_k\,.
$$
Let $a_\delta$
and $b_\delta$ be integers
such that
$a_\delta\delta\leq s_\delta-\delta^{\frac{1}{3}}<(a_\delta+1)\delta$
and
$b_\delta=a_\delta+[\delta^{-\frac{1}{6}}][\delta^{-\frac{1}{2}}]$,
where $[\cdot]$ denotes the integer part. By construction we have
that $a_\delta\delta\to t-$ and $b_\delta\delta\to t-$.
{}From the
estimate in Lemma~\ref{estim} between
$b_\delta$ and
$a_\delta$ we obtain
$$
\sum_{h=a_\delta+1}^{b_\delta} \|u_h^\delta-u_{h-1}^\delta\|^2 \le
M\,.
$$
Then
we divide the above sum
into
$[\delta^{-\frac{1}{6}}]$ groups of $[\delta^{-\frac{1}{2}}]$
consecutive
terms, and we find that the sum of one of these groups must be less
than or equal to $M/ [\delta^{-\frac{1}{6}}]$.
Therefore there exist two integers
$h_\delta$ and
$k_\delta$ such that
$a_\delta< h_\delta< k_\delta\le b_\delta$,
$k_\delta-h_\delta=[\delta^{-\frac{1}{2}}]$, and
$$
\sum_{h=h_\delta}^{k_\delta}\|u_h^\delta-u_{h-1}^\delta\|^2\le
\frac{M}{ [\delta^{-\frac{1}{6}}]}\,.
$$
It is then obvious that $k_\delta-h_\delta\to\infty$,
$h_\delta\delta$ and $k_\delta\delta$ converge to $t-$, and
(\ref{somma}) is satisfied.

Let us fix $s\in\Theta$ with $s<t$. Then
$$
K^\delta(s)\subset K^\delta(h_\delta\delta)\subset
K^\delta(k_\delta\delta)\subset K^\delta(s_\delta)\,,
$$
being $s<h_\delta\delta<k_\delta\delta\le s_\delta$ for $\delta $
small enough.
By compactness (Theorem~\ref{compactness}) we may assume that
$K^\delta(h_\delta\delta)\to K'$ and $K^\delta(k_\delta\delta)\to
K''$ in the Hausdorff metric.
Since
$K^\delta(s)\to K(s)$ and $K^\delta(s_\delta)\to K(t)$ in the Hausdorff
metric, we have
$$
K(s)\subset K'\subset K''\subset K(t)\,.
$$
Passing to the limit as $s\to t-$ we get
$$
K(t)=K(t-)\subset K'\subset K''\subset K(t)\,,
$$
which implies  that
$K^\delta(h_\delta\delta)$ and $K^\delta(k_\delta\delta)$
converge to $K(t)$ in the Hausdorff metric.
\end{proof}

\begin{lemma}\label{udeltat} For every $t\in\Theta$ we have
$u^\delta(t)-u^\delta(t-\delta)\to0$ strongly in $L^2(\Om)$ and $\nabla
u^\delta(t)\to\nabla u(t)$ strongly in $L^2(\Om;\R^2)$. Moreover
there exists a solution $u^*(t)$ of problem (\ref{egk}) corresponding
to $\phi=g(t)$ and $K=K(t)$ (possibly different from $u(t)$) such that
a subsequence of $u^\delta(t)$ (possibly depending on $t$) converges
to $u^*(t)$ strongly in $L^2(\Om)$.
\end{lemma}
\begin{proof}
Let $t\in\Theta$. By Lemma~\ref{l2} there exists a sequence
$t_\delta\to t-$  such that 
$K^\delta(t_\delta)\to K(t)$ in the Hausdorff metric
and
$u^\delta(t_\delta)-u^\delta(t_\delta-\delta)\to0$ strongly in~$L^2(\Om)$.
By the same argument,
we can
also construct $t'_\delta=l_\delta\delta$, with $l_\delta$ integer, such
that $t'_\delta\to t+$,
$K^\delta(t'_\delta)$ converge to
$K(t+)=K(t)$ in the Hausdorff metric, and
$u^\delta(t'_\delta)-u^\delta(t'_\delta-\delta)\to0$ strongly in~$L^2(\Om)$.

By Remark~\ref{linfty} the sequence
$u^\delta(t_\delta-\delta)$ is bounded in
$L^\infty(\Om)$, and by construction
$u^\delta(t_\delta)$ is the solution of the minimum problem 
(\ref{elambda}) which defines $\E_\lambda 
(g(t_\delta),K^\delta(t_\delta), u^\delta(t_\delta-\delta))$. 
Therefore Corollary~\ref{CC} implies that
$\nabla u^\delta(t_\delta)\to \nabla u(t)$ strongly 
in $L^2(\Om;\R^2)$. 
The same argument shows that
$\nabla u^\delta(t'_\delta)\to\nabla u(t)$ strongly in
$L^2(\Om;\R^2)$.

Then the estimate in Lemma~\ref{discr}
between
$h_\delta$ and $l_\delta$ gives
\begin{eqnarray*}
&\|\nabla u^\delta(t'_\delta)\|^2+
\huno(K^\delta(t'_\delta))+
\lambda\displaystyle
\sum_{h=h_\delta+1}^{l_\delta}\|u_h^\delta-u_{h-1}^\delta\|^2\le
\\ &\le \|\nabla u^\delta(t_\delta)\|^2+
\huno(K^\delta(t_\delta))+ 2 \displaystyle\int_{t_\delta}^{t'_\delta}
(\nabla u^\delta(\tau)|\nabla\gdot(\tau))d\tau+\rho(\delta)\,.
\end{eqnarray*}
Passing now to the limit as $\delta\to0$ we get
\begin{equation}\label{852002}
\sum_{h=h_\delta+1}^{l_\delta}\|u_h^\delta-u_{h-1}^\delta\|^2\to0\,.
\end{equation}
Let $i_\delta$ be the integer such that $i_\delta\delta\le
t<(i_\delta+1)\delta$. As $h_\delta<i_\delta\le l_\delta$, by
(\ref{852002}) we obtain that
$u^\delta(t)-u^\delta(t-\delta)= u_{i_\delta}^\delta-u_{i_\delta-1}^\delta
\to0$ strongly in $L^2(\Om)$.

Since $K^\delta(t)\to K(t)$ in the Hausforff metric and $u^\delta(t)$ 
is the solution of the minimum problem 
(\ref{elambda}) which defines
$\E_\lambda (g^\delta(t),K^\delta(t), u^\delta({t-\delta}))$,
from Corollary~\ref{CC} we obtain
that $\nabla u^\delta(t)\to
\nabla u(t)$ strongly in $L^2(\Om;\R^2)$ and that
 there exists a solution $u^*(t)$ of problem
(\ref{egk}) corresponding
to $\phi=g(t)$ and $K=K(t)$ (possibly different from $u(t)$)
such that 
a subsequence of $u^\delta(t)$ (possibly depending on $t$) converges
to $u^*(t)$ strongly in $L^2(\Om)$. 
\end{proof}

We show now that the increasing left-continuous function
$K\colon[0,T]\to\K$ satisfies
all conditions of Theorem~\ref{kt}. The following lemma proves 
condition~(b).

\begin{lemma}\label{condc}For every $t\in(0,T]$
we have
\begin{equation}\label{pt}
\E_\lambda(g(t),K(t),u(t))\leq \E_\lambda(g(t),K, u(t))\quad\forall\,
K\in\Kmf\,,\ K\supset K(t)\,.
\end{equation}
\end{lemma}
\begin{proof}
Let us consider first the case $t\in\Theta$. By Lemma \ref{udeltat},
$\nabla  u^\delta(t)\to\nabla u(t)$ strongly in
$L^2(\Om;\R^2)$ and,
passing  to a subsequence (which may depend on
$t$), we  may assume that $u^\delta(t)\to u^*(t)$ strongly in
$L^2(\Om)$, for some  solution $u^*(t)$ of the minimum problem
(\ref{egk}) corresponding to
$\phi=g(t)$ and $K=K(t)$. Then $\nabla u^*(t)=\nabla u(t)$ a.e.\ on 
$\Om$ and
$u^*(t)=u(t)$ a.e.\ on the connected
components of $\Omk(t)$ whose boundaries meet
$\partial_D\Omk(t)$, while on the
other connected components $u^*(t)$ and $u(t)$ are constant
(Remark~\ref{arm}). Moreover,
$\E_\lambda(g(t),K(t),u^*(t))=\E_\lambda(g(t),K(t),u(t))=\E(g(t),K(t))$.
By Lemma~\ref{udeltat} we have that $u^\delta({t-\delta})\to u^*(t)$
strongly in $L^2(\Om)$, and by Remark~\ref{linfty} 
$u^\delta({t-\delta})$ is bounded in $L^\infty(\Om)$.

Let $K\in\Kmf$ with $K\supset K(t)$.
Since $K^\delta(t)$ converges to $K(t)$ in the Hausdorff
metric as
$\delta\to0$, by \cite[Lemma 3.5]{DM-T} there exists a sequence
$K^\delta$ in $\Kmf$,
converging to $K$ in the Hausdorff metric, such that
$K^\delta\supset K^\delta(t)$ and
$\huno(K^\delta\setmeno K^\delta(t))\to
\huno(K\setmeno K(t))$
as $\delta\to0$.
By Lemma \ref{estim} this implies that $\huno(K^\delta)$ is
bounded as
$\delta\to0$.
Let $u^\delta$ and $u$ be the solutions of the minimum problems
(\ref{elambda}) which define
$\E_\lambda(g^\delta(t),K^\delta,u^\delta(t-\delta))$
and
$\E_\lambda(g(t),K,u^*(t))$, respectively. By Theorem~\ref{C}
$\nabla u^\delta\to\nabla u$ strongly in $L^2(\Om;\R^2)$.

The minimality of $K^\delta(t)$ expressed by (\ref{Kbul}) in Lemma~\ref{3.5}
gives
$$
\E_\lambda(g^\delta(t),K^\delta(t),
u^\delta(t-\delta))\leq
\E_\lambda(g^\delta(t),K^\delta,u^\delta(t-\delta))\,,
$$ which implies
$$
\|\nabla u^\delta(t)\|^2+
\lambda\,\|u^\delta(t)-u^\delta(t-\delta)\|^2\leq
\|\nabla u^\delta\|^2+\huno(K^\delta\setmeno K^\delta(t))
+\lambda\,\|u^\delta-u^\delta(t-\delta)\|^2\,.
$$
Passing now to the limit as $\delta\to0$
and using Lemma~\ref{udeltat} we get $\|\nabla u(t)\|^2\leq
\|\nabla u\|^2+
   \huno({K\setmeno K(t)})+\lambda\,\|u-u^*(t)\|^2$.
Adding $\huno(K(t))$ to both sides we obtain
$$
\E_\lambda(g(t),K(t),u(t))\le \E_\lambda(g(t),K,u^*(t))\,.
$$
As each connected component of $\Omk$ is contained in a connected
component of $\Omk(t)$, by Remark~\ref{costanti} we have
$\E_\lambda(g(t),K,u^*(t))=\E_\lambda(g(t),K,u(t))$. Therefore the previous
inequality gives (\ref{pt}) for
$t\in \Theta$.

Let us consider now the general case $t\in (0,T]$, which is
obtained by approximation.
We fix $t\in (0,T]$ and a compact set $K\in \Kmf$ with $K\supset
K(t)$. Let $t_k\to t-$, with $t_k\in\Theta$. Then $\nabla
u(t_k)\to\nabla u(t)$ strongly in $L^2(\Om;\R^2)$ by (\ref{ut-}). 
Arguing as in the proof of Theorem~\ref{C}, we may assume,  passing
to a subsequence,  that $u(t_k)\to u^*$
strongly in $L^2(\Om)$, for some  solution $u^*$ of the minimum
problem (\ref{egk}) corresponding to
$\phi=g(t)$ and $K=K(t)$. Then $\nabla u^*=\nabla u(t)$ a.e.\ on $\Om$ and
$u^*=u(t)$ a.e.\ on
the connected
components of $\Omk(t)$ whose boundaries meet
$\partial_D\Omk(t)$, while $u^*$ and $u(t)$ are constant on the other
connected components. Moreover
we have that
$\E_\lambda(g(t),K(t),u^*)=\E_\lambda(g(t),K(t),u(t))=
\E(g(t),K(t))$.
By the first part of the proof
$\E_\lambda(g(t_k),K(t_k),u(t_k))\le\E_\lambda(g(t_k),K,u(t_k))$.
Passing now to the limit as $k\to\infty$ thanks to Theorem~\ref{C} we get
\begin{eqnarray*}
&\displaystyle\E_\lambda(g(t),K(t),u(t))=\E_\lambda(g(t),K(t),u^*)
\le \liminf_{k\to\infty}\E_\lambda(g(t_k),K(t_k),u(t_k))\le\\
&\displaystyle\le\lim_{k\to\infty}\E_\lambda(g(t_k),K,u(t_k))=
\E_\lambda(g(t),K,u^*)=
\E_\lambda(g(t),K,u(t))\,,
\end{eqnarray*}
where the last equality follows from Remark~\ref{costanti}.
\end{proof}

The following lemma proves condition (c) of
Theorem~\ref{kt}.

\begin{lemma}\label{ineq}
For every $s,\, t$ with $0\leq s<t\leq T$
\begin{equation}\label{diseg}
\|\nabla u(t)\|^2+\huno(K(t))
\leq\|\nabla u(s)\|^2+\huno(K(s))+
2\displaystyle\int_s^t(\nabla u(\tau)|\nabla \gdot(\tau))d\tau\,.
\end{equation}
\end{lemma}
\begin{proof}
Let us fix $s,t\in \Theta$ with $0\leq s<t\leq T$.
Given $\delta>0$, let $i$ and $j$
be the
integers such that $t_i^\delta\leq s<t_{i+1}^\delta$ and
$t_j^\delta\leq t<t_{j+1}^\delta$.
Let us define $s_\delta:=t_i^\delta$ and $t_\delta:=t_j^\delta$.
Applying Lemma~\ref{discr} we obtain
\begin{equation}\label{ediscr}
\|\nabla u^\delta(t)\|^2+\huno(K^\delta(t)\setmeno K^\delta(s))
\leq\|\nabla u^\delta(s)\|^2+
\displaystyle2\int_{s_{\delta}}^{t_{\delta}}\!\!\!(\nabla u^\delta(\tau)|\nabla
\gdot(\tau))\,d\tau+\rho(\delta)\,,
\end{equation}
with $\rho(\delta)$ converging to zero as $\delta\to0$.
By Lemma~\ref{udeltat} for every $\tau\in\Theta$ we have $\nabla u^\delta(\tau)\to
\nabla u(\tau)$ strongly in
$L^2(\Om,\R^2)$ as $\delta\to0$, and by Lemma \ref{estim} we have
$\|\nabla
u^\delta(\tau)\|\le M$ for every $\tau\in[0,T]$. By
\cite[Corollary~3.4]{DM-T} we get
$$
\huno(K(t)\setmeno K(s))\leq\liminf_{\delta\to0}\,
\huno(K^\delta(t)\setmeno K^\delta(s))\,.
$$
Passing now to the limit in (\ref{ediscr}) as $\delta\to0$ we
obtain (\ref{diseg}) for every
$s, t\in\Theta$ with $0\leq s<t\leq T$.

In the general case we consider two sequences
$s_k\to s-$ and $t_k\to t-$ with $s_k, t_k\in\Theta$. Then
$K(s_k)\to K(s)$ and $K(t_k)\to K(t)$ in the Hausdorff metric, while
$\nabla u(s_k)\to\nabla u(s)$ and $\nabla u(t_k)\to\nabla u(t)$
strongly in $L^2(\Om;\R^2)$ by (\ref{ut-}).
By the first part of the proof we have that
\begin{equation}\label{disegk}
\|\nabla u(t_k)\|^2+\huno(K(t_k)\setmeno K(s_k))
\leq\|\nabla u(s_k)\|^2+
2\displaystyle\int_{s_k}^{t_k}(\nabla u(\tau)|\nabla
\gdot(\tau))\,d\tau\,.
\end{equation}
Passing now to the limit in (\ref{disegk}) as $k\to\infty$ and
using again
\cite[Corollary~3.4]{DM-T} we obtain~(\ref{diseg}).
\end{proof}

The following lemma proves condition (d) of Theorem~\ref{kt}.

\begin{lemma}\label{condb} We have
$
(u(0+),K(0+))\in{\mathcal R}^0(u_0,K_0)\,.
$
\end{lemma}
\begin{proof} By the definition of $u(0+)$ it follows 
that condition (a) in Definition~\ref{r0} is satisfied. We now take 
$(v_i^\delta, H_i^\delta):=(u_i^\delta,K_i^\delta)$.
By the argument used in the proof of Lemma
\ref{udeltat}, we can construct a sequence of integers 
$l_\delta\to\infty$ such that $l_\delta\delta\to 0+$,
$K^\delta(l_\delta\delta)$
converges to
$K(0+)$ in the Hausdorff metric,
and
$\nabla u^\delta(l_\delta\delta)\to\nabla u(0+)$ strongly in
$L^2(\Om;\R^2)$ as $\delta\to0$.
This proves that $(u(0+),K(0+)$ satisfies condition (b) in 
Definition~\ref{r0}.
\end{proof}

The following lemma proves condition (e) of Theorem~\ref{kt}.

\begin{lemma}\label{condd} For $0<t<T$ we have
$(u(t+),K(t+))\in{\mathcal R}^t(u(t),K(t))$.
\end{lemma}
\begin{proof} Fix $0<t<T$. By the definition of $u(t+)$ it follows 
that condition (a) in Definition~\ref{rt} is satisfied. We now take 
$(v_i^\delta, H_i^\delta):=(u_i^\delta,K_i^\delta)$.
Let $h_\delta$ and  $k_\delta$ be
the sequences of integers given by Lemma \ref{l2}.
Since $\sum_{h=h_\delta}^{k_\delta}
\|u_h^\delta-u_{h-1}^\delta\|^2\to0$
as $\delta\to0$, we have
$u_{\sigma_\delta}^\delta-
u_{{\sigma_\delta}-1}^\delta\to0$ strongly in $L^2(\Om)$ 
for every sequence $\sigma_\delta$ of integers between
$h_\delta$ and
$k_\delta$.
Since both $K^\delta(h_\delta\delta)$ and
$K^\delta(k_\delta\delta)$ converge to
$K(t)$ in the Hausdorff metric,
$K^\delta({\sigma_\delta}\delta)$
converges to
$K(t)$ in the Hausdorff metric. Therefore
$\nabla u^\delta({\sigma_\delta}\delta)\to\nabla u(t)$ strongly in 
$L^2(\Om;\R^2)$ by Corollary~\ref{CC}.
This shows that  condition $({\rm b}_1)$  in
Definition~\ref{rt}  is satisfied.

By the same
argument as in the proof of Lemma \ref{udeltat},  we can
construct a sequence of integers $l_\delta\to\infty$
such that $l_\delta\delta\to t+$, $l_\delta-k_\delta\to\infty$, 
$K^\delta(l_\delta\delta)$
converges to $K(t+)$ in the Hausdorff metric, and
$\nabla  u^\delta(l_\delta\delta)\to\nabla u(t+)$ strongly in $L^2(\Om;\R^2)$.
This shows that  condition $({\rm b}_2)$ in
Definition~\ref{rt} is satisfied.
\end{proof}

\end{section}

\begin{section}{Example}\label{esempio}

In this section we consider in detail the particular case when no initial crack
is present, i.e., $K_0=\emptyset$. We prove that, if $\Om$ and $g(t)$
are sufficiently regular, no crack will appear in our model, provided
$\lambda $ is large enough. More precisely, under these conditions we
prove that $K(t)=\emptyset$ is the unique function which satisfies
  conditions (a)--(e) of Theorem~\ref{kt}.

\begin{theorem}\label{ex}
Assume that $\partial\Om$ is of class $C^2$,
$\partial_D\Om=\partial\Om$,
   and
$g\in AC([0,T];H^1(\Om))\cap L^\infty([0,T];C^{1,\alpha}(\overline\Om))\cap
C^0([0,T];C^0(\overline\Om))$ for some $0<\alpha<1$.
If $K_0=\emptyset$ and
$\lambda$ is larger than the constant $\lambda_0$ given by (\ref{lambda0}), 
then
$K(t)=\emptyset$ is the unique function $K\colon[0,T]\to\Kmf$ which
satisfies conditions (a)--(e) of Theorem~\ref{kt}. Moreover,
\begin{equation}\label{Rempty}
{\mathcal R}^t(u(t),\emptyset)=\{(u(t),\emptyset)\}
\end{equation}
for every $t\in[0,T)$.
\end{theorem}

To prove Theorem~\ref{ex} we need some estimates on the solutions of
the Dirichlet problems
\begin{equation}\label{P1}
\left\{\begin{array}{ll} \Delta u=f&
\hbox{in }\Om\,,\\
\salt
u=\phi &\hbox{on } \partial\Om\,,
\end{array}\right.
\end{equation}
and
\begin{equation}\label{P3}
\left\{\begin{array}{ll} \Delta v=\lambda(v-w) &
\hbox{in }\Om\,,\\
\salt
v=\psi &\hbox{on } \partial\Om\,.
\end{array}\right.
\end{equation}
If $\partial\Om\in C^{1,\alpha}$ and  $\phi\in
C^{1,\alpha}(\overline\Om)$, for $0<\alpha<1$, and $f\in
L^\infty(\Om)$, then the solution $u$ of (\ref{P1}) belongs to
$C^{1,\alpha}(\overline\Om)$ (see, e.g., \cite[Corollary 8.35]{G-T})
and there exists a constant $C$, 
independent  of $f$ and $\phi$, such that 
\begin{equation}\label{8.35}
\|\nabla u\|_{\infty}\le
C\,(\|f\|_{\infty}+\|\nabla\phi\|_{0,\alpha})\,,
\end{equation}
where 
$\|\cdot\|_{0,\alpha}$ denotes the norm in
$C^{0,\alpha}(\overline\Om;\R^2)$ and $\|\cdot\|_\infty$ denotes the
norm in $L^\infty(\Om)$ or in $L^\infty(\Om;\R^2)$, according to the
context.

If $w\in L^2(\Om)$ and $\psi\in H^1(\Om)\cap L^\infty(\Om)$, then the
solution $v$ of (\ref{P3}) belongs to $H^1(\Om)\cap L^\infty(\Om)$
and
\begin{equation}\label{8.16}
     \|v\|_{\infty}\le
     \|\psi\|_{\infty}+C_\lambda\,\|w\|\,,
     \end{equation}
where the constant $C_\lambda$ depends on $\lambda$, but  not
 on $w$ and $\psi$ (see, e.g., \cite[Theorem 8.16]{G-T}).
\begin{proof}[Proof of Theorem \ref{ex}.]
We begin by proving that $K(t)=\emptyset$ satisfies condition (b).
Since every $K\in\Kmf$ can be approximated in the Hausdorff metric
by  a
sequence of compact sets contained in $\Om$, with convergence of the
lenghts, taking Theorem~\ref{C} into account 
it is enough to prove that for
every $0<t\le T$ we have
\begin{equation}\label{vuoto}
\E_\lambda(g(t),\emptyset,u(t))\le\E_\lambda(g(t),K,u(t))
  \end{equation}
  for every compact set $K\subset\Om$.
  
  To this end we use the calibration constructed in \cite[Section
  5.3]{A-B-DM}. In  that  section the Neumann condition on $\partial\Om$
  is  used only to obtain that $\phi^x(x,t)\,\nu=0$ for $x\in\partial\Om$,
  which is not needed in our case, where we prescribe a Dirichlet
  boundary condition on $\partial\Om$ (see \cite[Theorem 3.3]{A-B-DM}).
This calibration can be constructed provided we are able to prove
inequality (5.12) of
\cite{A-B-DM}, which in our case reduces to
\begin{equation}\label{128}
2^7\,\|\nabla u(t)\|_{\infty}^4<\lambda\,.
\end{equation}
By (\ref{8.35}) there exists a constant $C$
such that
$$
\|\nabla u(t)\|_{\infty}\le C\, G_{\alpha}\,,
$$
where
$$
G_{\alpha}:=\sup_{t\in[0,T]}\|\nabla g(t)\|_{0,\alpha}\,.
$$
Therefore (\ref{128}) is satisfied if
\begin{equation}\label{lambda0}
\lambda>\lambda_0:=2^7\, C^4\, G_{\alpha}^4\,,
\end{equation}
and in this case the calibration constructed in 
\cite[Section 5.3]{A-B-DM} proves (\ref{vuoto})
for every compact set $K\subset\Om$, which implies 
condition (b) of Theorem~\ref{kt}.

Let us prove now that
${\mathcal R}^0(u(0),\emptyset)=\{(u(0),\emptyset)\}$. As  
${\mathcal R}^0(u(0),\emptyset)\neq\emptyset$ (see, e.g., 
Lemma~\ref{condb}), it is enough to show that 
${\mathcal R}^0(u(0),\emptyset)\subset\{(u(0),\emptyset)\}$.

Let $(v,H)\in{\mathcal R}^0(u(0),\emptyset)$, let $\delta_n$, 
$l_n$, $v_i^{\delta_n}$, and $H_i^{\delta_n}$ be the sequences, the 
functions, and the sets which appear in  condition (b) of Definition~\ref{r0},
and let
\begin{equation}\label{cont}
    \e_n:= \sup_{t\in[\delta_n,T]}
    \|g(t)-g(t-\delta_n)\|_\infty\,.
\end{equation}
As $g\in C^0([0,T], C^0(\overline\Om))$,  the sequence $\e_n$ tends 
to~$0$.
Starting from $w_0^{\delta_n}:=u(0)$, we consider  also the  sequence
$w_i^{\delta_n}$, $0<i\le l_n$,
of  the solutions of the
Dirichlet problems
\begin{equation}\label{Pi}
\left\{\begin{array}{ll} \Delta 
w_i^{\delta_n}=\lambda(w_i^{\delta_n}-
w_{i-1}^{\delta_n}) &
\hbox{in }\Om\,,\\
\salt
   w_i^{\delta_n}=g_i^{\delta_n} &\hbox{on } \partial\Om\,.
\end{array}\right.
\end{equation}

By using the calibration constructed in \cite[Section 5.3]{A-B-DM}, we will
prove by induction on $i$ that
$H_i^{\delta_n}=\emptyset$ and $v_i^{\delta_n}= w_i^{\delta_n}$ for 
$\lambda>\lambda_0$ and for $n$ large enough.
This calibration can be constructed provided we are able to prove
inequality (5.12) of
\cite{A-B-DM}, which in this case reads
\begin{equation}\label{128a}
\|\nabla w_i^{\delta_n}\|_\infty(\sqrt\lambda\,\|w_i^{\delta_n}-
w_{i-1}^{\delta_n}\|_\infty+\sqrt2\,\|\nabla
w_i^{\delta_n}\|_\infty)<\frac{1}{8}\sqrt\lambda\,.
\end{equation}
To obtain (\ref{128a}) we prove by induction on $i$ that
\begin{equation}\label{Delta}
     \|\Delta w_i^{\delta_n}\|_\infty\le\lambda\,\e_n\,.
   \end{equation}
This inequality is true for $i=0$ since $w_0^{\delta_n}=u(0)$, which is
harmonic. Assume that (\ref{Delta}) holds for $i-1$. Then
$w:=w_{i-1}^{\delta_n}+\e_n$ is a super-solution of  
problem (\ref{Pi}), 
 in the sense that
$$
\left\{\begin{array}{ll} \Delta w\le\lambda(w-w_{i-1}^{\delta_n}) &
\hbox{in }\Om\,,\\
\salt
w\ge g_i^{\delta_n}&\hbox{on } \partial\Om\,.
\end{array}\right.
$$
Indeed,  $\Delta w=\Delta w_{i-1}^{\delta_n}\le\lambda\,\e_n=
\lambda(w-w_{i-1}^{\delta_n}) $ in 
$\Om$ by the
inductive hypothesis, while
$w=g_{i-1}^{\delta_n}+\e_n
\ge g_i^{\delta_n}$ on
$\partial\Om$ by (\ref{cont}).
Therefore $w_i^{\delta_n}\le w_{i-1}^{\delta_n}+\e_n$ in $\Om$.
Similarly, $w_{i-1}^{\delta_n}-\e_n$ is a sub-solution of 
(\ref{Pi}); this implies that $w_i^{\delta_n}\ge 
w_{i-1}^{\delta_n}-\e_n$, which, together with the 
previous inequality gives
\begin{equation}\label{ii+1}
     \|w_i^{\delta_n}-w_{i-1}^{\delta_n}\|_\infty\le\e_n\,.
\end{equation}
By (\ref{Pi})  we have
$\|\Delta
w_i^{\delta_n}\|_\infty\le \lambda\,\e_n$, concluding the proof of
(\ref{Delta}).

{}From (\ref{8.35}) and (\ref{Delta}) 
we obtain
\begin{equation}\label{nablaui}
\|\nabla w_i^{\delta_n}\|_\infty\le
C\,(\lambda\,\e_n +G_{\alpha})\,.
\end{equation}
By (\ref{ii+1}) and (\ref{nablaui}) for $n$ large 
enough we have
\begin{eqnarray*}
 &\displaystyle \|\nabla w_i^{\delta_n}\|_\infty(\sqrt\lambda\,\|w_i^{\delta_n}-
w_{i-1}^{\delta_n}\|_\infty+\sqrt2\,
\|\nabla w_i^{\delta_n}\|_\infty)
\le
\\
& \displaystyle \le C\,(\lambda\,\e_n +
G_{\alpha}) \, (\sqrt\lambda\,\e_n + \sqrt2\,
C\,(\lambda\,\e_n +G_{\alpha})) < \frac{1}{8}\sqrt\lambda\,,
\end{eqnarray*}
where the last inequality follows from  (\ref{lambda0}) and from the 
fact that ${\e_n\to0}$. This proves (\ref{128a})  for $n$ large 
enough.

Therefore, using the calibration constructed in
\cite[Section 5.3]{A-B-DM}, we can prove that
$(w_i^{\delta_n},\emptyset)$ is
the unique minimizer of the functional
$$
F_{i-1}^{\delta_n}(u,K):=\|\nabla u\|^2+\huno(K)+
\lambda\,\|u-w_{i-1}^{\delta_n}\|^2
$$
among all pairs $(u,K)$ with $K\in \Kmfbul$, $u\in H^1(\Omk)$, 
$u=g_i^{\delta_n}$ on $\partial\Omk$.

As $v_0^{\delta_n}=u(0)=w_0^{\delta_n}$ and 
$H_0^{\delta_n}=K(0)=\emptyset$, both $(v_1^{\delta_n},H_1^{\delta_n})$ and
$(w_1^{\delta_n},\emptyset)$ minimize $F_0^{\delta_n}$ with the same 
Dirichlet condition
$g_1^{\delta_n}$, hence the  uniqueness result gives $v_1^{\delta_n}=
w_1^{\delta_n}$ and
$H_1^{\delta_n}=\emptyset$. In the same way, by induction we prove that
$v_i^{\delta_n}=w_i^{\delta_n}$ and $H_i^{\delta_n}=\emptyset$ for 
every~$i$.

By condition $({\rm b}_2)$ of Definition~\ref{r0} we have $H=\emptyset$ and
$\nabla w_{l_n}^{\delta_n}\to \nabla v$ strongly in 
$L^2(\Om;\R^2)$. As $w_{l_n}^{\delta_n}-w_{l_n-1}^{\delta_n}\to0$ 
strongly in $L^2(\Om)$ by (\ref{ii+1}) and
$g_{l_n}^{\delta_n}=g(l_n\delta_n)\to g(0)$ 
strongly in $H^1(\Om)$,
the continuous dependence of the solutions of (\ref{Pi})
on the data implies
that $w_{l_n}^{\delta_n}$ converges to $u(0)$ strongly in $H^1(\Om)$.
This shows that 
$v=u(0)$ and concludes the proof of the inclusion
${\mathcal R}^0(u(0),\emptyset)\subset\{(u(0),\emptyset)\}$.

Let us prove now that
${\mathcal R}^t(u(t),\emptyset)=\{(u(t),\emptyset)\}$ for every $0<t<T$.
As  
${\mathcal R}^t(u(t),\emptyset)\neq\emptyset$ (see, e.g., 
Lemma~\ref{condd}), it is enough to show that 
${\mathcal R}^t(u(t),\emptyset)\subset\{(u(t),\emptyset)\}$.
Let $(v,H)\in{\mathcal R}^t(u(t),\emptyset)$ and let $\delta_n$, 
$h_n$, $k_n$,  $l_n$, $v_i^{\delta_n}$, and $H_i^{\delta_n}$ be the 
sequences, the functions and the sets which appear in condition (b) of 
Definition~\ref{rt}.
As $\emptyset$ is isolated in the Hausdorff metric, by $({\rm b}_1)$ we may 
assume that $H_{k_n}^{\delta_n}=\emptyset$, and by the 
monotonicity of $H_i^{\delta_n}$ we deduce that $H_i^{\delta_n}=\emptyset$ 
for $0\le i\le k_n$.
It follows that $v_i^{\delta_n}$ belongs to  $H^1(\Om)$ and
solves
  the Dirichlet problem
\begin{equation}\label{Pi*}
\left\{\begin{array}{ll} \Delta v_i^{\delta_n}=
\lambda(v_i^{\delta_n}-v_{i-1}^{\delta_n}) &
\hbox{in }\Om\,,\\
\salt
    v_i^{\delta_n}=g_i^{\delta_n} &\hbox{on } \partial\Om\,,
\end{array}\right.
\end{equation}
for $1\le i\le k_n$ .

In order to prove that $H_i^{\delta_n}=\emptyset$ for $ k_n<i\le l_n$
we 
apply the calibration method as in the case $t=0$. 
We define $w_{k_n}^{\delta_n}:=v_{k_n}^{\delta_n}$ and 
we consider the sequence $w_i^{\delta_n}$, $k_n<i\le l_n$, 
defined inductively by the solutions of 
 (\ref{Pi}).
We can construct a calibration for $w_i^{\delta_n}$ provided (\ref{128a})
holds. As before,
it is enough to show  that $\Delta w_i^{\delta_n}\to 0$  
in $L^\infty(\Om)$ as $n\to \infty$, uniformly for $k_n\le i\le l_n$. 
The inductive argument used to prove (\ref{Delta}) shows that
$$
\|\Delta w_{i}^{\delta_n}\|_\infty\le 
\max\{\lambda\,\e_n, 
\|\Delta w_{k_n}^{\delta_n}\|_\infty\}=
\max\{\lambda\,\e_n, 
\|\Delta v_{k_n}^{\delta_n}\|_\infty\}
$$
for $k_n\le i\le l_n$.
Therefore, in order to obtain (\ref{128a}), it is enough to show that
$\Delta v_{k_n}^{\delta_n}\to 0$  in $L^\infty(\Om)$.

By
(\ref{Pi*}) we have
$\Delta 
v_{k_n}^{\delta_n}=
\lambda(v_{k_n}^{\delta_n}-v_{k_n-1}^{\delta_n})$. To
estimate $\|v_{k_n}^{\delta_n}-v_{k_n-1}^{\delta_n}\|_\infty$ 
we note that the
difference  satisfies
\begin{equation}\label{Pdiff}
\left\{\begin{array}{ll}
\Delta(v_{k_n}^{\delta_n}-v_{k_n-1}^{\delta_n})=
\lambda((v_{k_n}^{\delta_n}-v_{k_n-1}^{\delta_n})-
(v_{k_n-1}^{\delta_n}-v_{k_n-2}^{\delta_n})) &
\hbox{in }\Om\,,\\
\salt
   v_{k_n}^{\delta_n}-v_{k_n-1}^{\delta_n}=
   g_{k_n}^{\delta_n}-g_{k_n-1}^{\delta_n} 
   &\hbox{on } \partial\Om\,.
\end{array}\right.
\end{equation}
As $g_{k_n-1}^{\delta_n}-g_{k_n-2}^{\delta_n} \to0$ strongly in 
$H^1(\Om)$, and 
$\nabla v_{k_n-1}^{\delta_n}-\nabla v_{k_n-2}^{\delta_n} \to 
0$ strongly in $L^2(\Om;\R^2)$ (by condition $({\rm b}_1)$ of 
Definition~\ref{rt}), using the  Poincar\'e inequality we conclude 
that $v_{k_n-1}^{\delta_n}-v_{k_n-2}^{\delta_n} \to 0$ strongly in 
$L^2(\Om)$. 

Since $g_{k_n}^{\delta_n}-g_{k_n-1}^{\delta_n} \to0$ in 
$L^\infty(\Om)$, estimate (\ref{8.16}) for (\ref{Pdiff}) implies that 
 $v_{k_n}^{\delta_n}-v_{k_n-1}^{\delta_n}\to0$ in $L^\infty(\Om)$. 
By (\ref{Pi*}) this implies that $\Delta v_{k_n}^{\delta_n}\to 0$ 
in $L^\infty(\Om)$.

Therefore, arguing  as in the case $t=0$, we can construct now a 
calibration for 
$w_i^{\delta_n}$, which shows that $H_i^{\delta_n}=\emptyset$ 
and $v_i^{\delta_n}=w_i^{\delta_n}$ for $k_n\le i\le l_n$, 
and leads to the conclusion of the proof of (\ref{Rempty}). 

So far we have proved that $K(t)=\emptyset$ satisfies conditions 
(a), (b), (d), and (e) of Theorem~\ref{kt}.  
As $\E(g(t), K(s))=\E(g(t),\emptyset)$, 
condition ($ \mskip-.7\thinmuskip{\rm h}'\mskip-.7\thinmuskip$) 
of Proposition~\ref{gtKs} is trivial. Condition 
(f) of Remark~\ref{3.8}  follows from the smooth dependence of the energy
on the boundary data. By Proposition~\ref{gtKs} conditions (f) and 
($ \mskip-.7\thinmuskip{\rm h}'\mskip-.7\thinmuskip$) imply 
condition (c)  of Theorem~\ref{kt}.

Let us prove now the uniqueness.
Let
$\tilde K\colon[0,T]\to\Kmf$ be another function which
satisfies conditions (a)--(e) of Theorem~\ref{kt}, and let
$\tilde u(t)$ be a solution of the minimum problem (\ref{egk})
corresponding to $\phi=g(t)$ and
$K=\tilde K(t)$. Assume by contradiction that there exists an
instant $t\in[0,T]$ such that $\tilde K(t)\neq\emptyset$
and let $t_0$ be the infimum of such instants. By the finite
intersection property we have $\tilde K(t_0+)\neq\emptyset$. We will show
  that properties (a), (d), and (e), together with (\ref{Rempty}),
  imply that $\tilde K(t_0+)=\emptyset$. This contradiction proves
  that $\tilde K(t)=\emptyset$ for every $t\in[0,T]$.

  If $t_0=0$, by properties (a) and  (d), and by (\ref{Rempty})
  we have
  $$
  (\tilde u(0+),\tilde K(0+))\in{\mathcal R}^0(\tilde u(0),\tilde K(0))=
  {\mathcal R}^0( u(0),\emptyset)=\{(u(0),\emptyset)\}\,,
  $$
hence $\tilde K(0+)=\emptyset$.

If $t_0>0$, we have $\tilde K(t)=\emptyset$ and $\tilde u(t)=u(t)$
for $0\le t<t_0$. Hence $\tilde K(t_0-)=\emptyset$ and $\tilde
u(t_0-)=u(t_0-)$. By property (e) and by (\ref{Rempty})
  we have
  $$
  (\tilde u(t_0+),\tilde K(t_0+))\in{\mathcal R}^{t_0}(\tilde
u(t_0-),\tilde K(t_0-))=
  {\mathcal R}^{t_0}( u(t_0-),\emptyset)=\{(u(t_0),\emptyset)\}\,,
  $$
hence $\tilde K(t_0+)=\emptyset$. This concludes the proof of the
uniqueness.
\end{proof}

\end{section}

\begin{section}{Behaviour Near the Tips}\label{tips}

 In this section, given  $g\in
AC([0,T];H^1(\Om))\cap L^\infty([0,T];L^\infty(\Om))$, we consider a 
function $K\colon[0,T]\to\Kmf$
which satisfies conditions (a)--(e) of Theorem~\ref{kt}, and  study
the behaviour of the solutions $u(t)$ near the ``tips'' of the sets
$K(t)$. Under some natural assumptions on the geometry of the sets $K(t)$, 
we shall see that $K(t)$
satisfies Griffith's criterion for crack growth.   

More precisely, let
$0\le t_0<t_1\le T$. Suppose that the following structure condition
is satisfied: there exists a finite family of simple arcs
${\Gamma}_i$,
$i=1,\,\ldots,p$, contained in $\Om$ and parametrized by arc length
by $C^2$  paths $\gamma_i\colon [\sigma_i^0, \sigma_i^1]\to \Om$,
such that, for $t_0<t<t_1$,
\begin{equation}\label{structure}
K(t)=K(t_0)\cup \bigcup_{i=1}^p {\Gamma}_i(\sigma_i(t))\,,
\end{equation}
where ${\Gamma}_i(\sigma):=
\{\gamma_i(\tau): \sigma_i^0\le \tau \le \sigma\}$ and
$\sigma_i\colon [t_0,t_1]\to [\sigma_i^0, \sigma_i^1]$ are
nondecreasing functions with $\sigma_i(t_0)=\sigma_i^0$ and
$\sigma_i^0<\sigma_i(t)<\sigma_i^1$ for $t_0<t<t_1$.
Assume also that the arcs
${\Gamma}_i$ are pairwise disjoint, and that
${\Gamma}_i\cap K(t_0)=\{\gamma_i(\sigma_i^0)\}$.
For $i=1,\,\ldots,p$ and
$\sigma_i^0<\sigma<\sigma_i^1$
let $\kappa_i(u,\sigma)$ be the stress intensity factor
defined by (8.2) in \cite{DM-T} with $\gamma=\gamma_i$ and $B$ equal to a
sufficiently small
ball centred at $\gamma_i(\sigma)$.

We are now in a position to state the main result of this section, 
which expresses Griffith's criterion in our model.

\begin{theorem}\label{Griffith}
Let $T>0$, $\lambda>0$, $m\ge1$,  and $g\in AC([0,T];H^1(\Om))\cap 
L^\infty([0,T];L^\infty(\Om))$. Let $K\colon[0,T]\to\Kmf$ be a function
which satisfies conditions (a)--(e) of Theorem~\ref{kt}, and  let
$u(t)$ be a solution of the minimum problem (\ref{egk})  defining
$\E(g(t),K(t))$. Given $0\le t_0<t_1\le T$,
assume that (\ref{structure}) is satisfied for $t_0<t<t_1$, and that
the arcs ${\Gamma}_i$ and the functions $\sigma_i$
satisfy all properties considered above.
Then
\begin{eqnarray}
&\dot\sigma_i(t)\ge 0\quad \hbox{ for a.e.\ } t\in (t_0,t_1)\,,
\label{sigmadot}\\
&1-\kappa_i(u(t),\sigma_i(t))^2\ge 0
\quad \hbox{ for every\ } t\in (t_0,t_1)\,,
\label{sif>}\\
&\big\{1- \kappa_i(u(t),\sigma_i(t))^2\big\}
\,\dot\sigma_i(t)= 0\quad \hbox{ for a.e.\ } t\in (t_0,t_1)\,,
\label{sif=}
\end{eqnarray}
for $i=1,\,\ldots,p$.
\end{theorem}

The proof of Theorem~\ref{Griffith} is obtained by adapting the proof 
of Theorem~8.4 of~\cite{DM-T}. We indicate here only the changes to be 
done. 

First of all, we need a localized version of the 
energies $\E$ and $\E_\lambda$. 
If $A$ is a bounded open set in $\R^2$ with Lipschitz boundary,
$K$ is a compact set in $\R^2$, 
$\phi\colon \partial A\setmeno K\to\R$ is a bounded function, and 
$w\in L^2(A)$, we define
\begin{eqnarray}
&\qquad\E(\phi,K, A):=\displaystyle\min_{v\in{\mathcal V}(\phi,K,A)}
\Big\{\int_{A\setminus K}|\nabla v|^2\,dx
+\huno(K\cap \overline A)\Big\} \,,\label{egkloc}\\
&\qquad\E_\lambda(\phi,K,A,w):=\displaystyle\min_{v\in{\mathcal
V}(\phi,K,A)}
\Big\{\int_{A\setminus K}|\nabla v|^2\,dx
+\huno(K\cap \overline A)+\lambda\int_{A\setminus K}|
v-w|^2\,dx\Big\} \,,\label{elambdaloc}
\end{eqnarray}
where
$$
{\mathcal V}(\phi,K,A):=\{v\in H^1(A\setmeno K):
v=\phi\quad\hbox{ on }
\partial A\setmeno K\}\,.
$$

Then we can prove the  following result for $\E_\lambda$, arguing  as in 
\cite[Lemma 8.5]{DM-T}. 

\begin{lemma}\label{cond3loc'} Let $m\ge1$, $\lambda>0$, let
$H\in\Kmf$ with
$h$ connected components,
let $\phi\in H^1(\Om)$, $w\in L^2(\Om)$, and let $u$ be the solution
of the minimum problem (\ref{elambda}) which
defines $\E_\lambda(\phi,H,w)$.
   Given  an open subset $A$ of $\Om$, with
Lipschitz boundary, such that $H\cap\overline A\neq\emptyset$, let
$q$ be the number of connected components of $H$ which meet
$\overline A$.
Assume that
\begin{equation}\label{minnl}
\E_\lambda(\phi,H,w)\leq \E_\lambda(\phi,K,w)\qquad\forall
\,K\in\Kmf,\,\ K\supset H\,.
\end{equation}
Then
\begin{equation}\label{minnlA}
\E_\lambda(u,H,A,w)\leq \E_\lambda(u,K,A,w)\qquad\forall
\,K\in\KA,\,\ K\supset H\cap\overline A\,.
\end{equation}
\end{lemma}

\begin{proof}[Proof of Theorem~\ref{Griffith}.]
We now consider in detail the changes needed in the proof of Theorem~8.4 of 
\cite{DM-T}. Inequality (8.12) must be 
replaced by 
\begin{equation}\label{ge0}
\frac{d}{d\sigma} \E_\lambda(u(t),{\Gamma}_i(\sigma),B_i,u(t))
\Big|\lower1.5ex\hbox{$\scriptstyle \sigma=\sigma_i(t)$}
\ge 0\,,
\end{equation}
which can be derived, arguing as in \cite{DM-T}, from  Lemma~\ref{cond3loc'}  
and from the minimality 
property (b) of Theorem~\ref{kt}. 

On the other hand, we can show that 
$$
\frac{d}{d\sigma} \E_\lambda(u(t),{\Gamma}_i(\sigma),B_i,u(t))
\Big|\lower1.5ex\hbox{$\scriptstyle 
\sigma=\sigma_i(t)$}=1-\kappa_i(u(t),\sigma_i(t))^2
$$
by adapting the proof of \cite[Theorem 6.4.1]{Gri2}.
This equality, together with (\ref{ge0}), proves (\ref{sif>}).

To obtain (\ref{sif=})  we  continue the proof of Theorem~8.4 of 
\cite{DM-T}, noting that 
 the  inequality in condition 
 ($ \mskip-.7\thinmuskip{\rm h}'\mskip-.7\thinmuskip$) 
 of Proposition~\ref{gtKs} is enough to 
 conclude the proof.
\end{proof}

\medskip

{\bf Acknowledgements.} The work of Gianni Dal Maso is part of the 
Project ``Calculus of Variations'',
supported by SISSA and by
the Italian Ministry of Education, University, and Research.
The work of Rodica Toader is part of the Project ``Methods and 
Problems in Real Analysis'',
supported by the University of Trieste and by
the Italian Ministry of Education, University, and Research.
\end{section}

\end{document}